\documentclass[preprint,12pt]{elsarticle}
\RequirePackage{etex}
\RequirePackage{tikz}
\RequirePackage{graphicx}
\RequirePackage{pgfplots}
\RequirePackage{epstopdf}
\usepgfplotslibrary{groupplots}
\pgfplotsset{compat=newest}
\usepackage{bigints}
\pgfplotsset{
    every axis plot/.append style = {font = \scriptsize}
  }

\RequirePackage{etoolbox}

\RequirePackage{ifthen}   

\RequirePackage{amsfonts}
\RequirePackage{amsthm}
\RequirePackage{amsmath}
\RequirePackage{mathtools,xparse}
\RequirePackage[utf8]{inputenc}
\RequirePackage{url}

\RequirePackage{marvosym}
\RequirePackage{bm}
\RequirePackage[font=small,labelfont=bf]{caption}
\RequirePackage[ruled]{algorithm2e}
\RequirePackage{comment}
\RequirePackage{mathrsfs}
\RequirePackage{stmaryrd}

\RequirePackage{hyperref}

\RequirePackage{float}

\usetikzlibrary{quotes,angles,positioning}

\RequirePackage{caption}
\RequirePackage{subcaption}

\RequirePackage{amssymb}
\usetikzlibrary{matrix}

\newtheorem{definition}{Definition}[section]

\newtheorem{assumption}{Assumption}

\newcommand{\rr}{\mathbb{R}}
\newcommand{\pp}{\mathbb{P}}

\newcommand{\Pidel}{\Pi^{\nabla}}
\newcommand{\PiO}{\Pi^{0}}
\newcommand{\PiI}{\Pi^{1}}
\newcommand{\Jacobian}{\mathbf{J}}

\newcommand{\RefGrad}{\hat{\nabla}}

\newcommand{\RefCoordinate}{\boldsymbol{\xi}}
\newcommand{\PhysCoordinate}{\mathbf{x}}

\newcommand{\RefDomain}{\hat{\Omega}}

\newcommand{\ALEVelocity}{\mathbf{w}}

\newcommand{\ALEMap}{\mathscr{A}}
\newcommand{\RefMesh}{\hat{\mathcal{T}}_h}

\newcommand{\norm}[1]{\left\| #1 \right\|}
\newcommand{\abs}[1]{\left| #1 \right|}
\renewcommand{\det}[1]{\text{det} \left( #1 \right)}

\newcommand{\dofs}{\text{dofs}}
\newcommand{\Ndofs}{N^{\dofs}}

\renewcommand{\b}{\mathbf{b}}
\newcommand{\RefE}{\hat{E}}

\journal{Computers \& Mathematics with Applications}

\begin{document}

\title{A High-order Arbitrary Lagrangian-Eulerian Virtual Element Method for Convection-Diffusion Problems}

\begin{frontmatter}

\author[inst1]{H. Wells \corref{cor1}}
\cortext[cor1]{Corresponding author}
\ead{h.wells.research@gmail.com}

\affiliation[inst1]{organization={School of Mathematical Sciences, University of Nottingham},
            addressline={University Park}, 
            city={Nottingham},
            postcode={NG7 2RD}, 
            country={United Kingdom}}

\begin{abstract}
    A virtual element discretisation of an Arbitrary Lagrangian-Eulerian method for two-dimensional convection-diffusion equations is proposed employing an isoparametric Virtual Element Method to achieve higher-order convergence rates on curved edged polygonal meshes. The proposed method is validated with numerical experiments in which optimal $H^1$ and $L^2$ convergence are observed. This method is then successfully applied to an existing moving mesh algorithm for implicit moving boundary problems in which higher-order convergence is achieved.
\end{abstract}

\begin{keyword}
Moving Mesh Method \sep Virtual Element Method \sep Arbitrary Lagrangian-Eulerian Schemes \sep Polygonal Meshes \sep Convection-diffusion Equations
\PACS 0000 \sep 1111
\MSC 0000 \sep 1111
\end{keyword}

\end{frontmatter}

\section{Introduction}\label{sec::introduction}

In Computational Fluid Dynamics (CFD), Arbitrary Lagrangian-Eulerian (ALE) schemes are pivotal for accurately simulating fluid flows, especially in scenarios involving large deformations or moving boundaries \cite{Hirt1974AnSpeeds,Donea1982AnInteractions,Donea2004}. By combining aspects of both Lagrangian and Eulerian methods, ALE schemes offer improved mesh adaptability and computational efficiency. Their versatility makes them suitable for a variety of CFD applications, from aerodynamics to fluid-structure interactions \cite{Takashi1992AnBody,Souli2000ALEProblems,Wesseling2001PrinciplesDynamics}.

Common challenges in ALE schemes, such as mesh tangling and distortion during substantial deformations, sustaining mesh quality over time, and ensuring numerical stability and accuracy, are well-documented \cite{Donea2004}. Effectively addressing these concerns is vital for the successful application of ALE schemes in intricate fluid simulations. The adoption of polygonal discretization techniques facilitates the use of more distorted and irregular elements while preserving numerical stability. Furthermore, polygonal discretization allows for efficient resolution of complex geometries with fewer computational elements, achieving comparable numerical precision to traditional triangular and quadrilateral meshes, for instance, through discontinuous Galerkin methods \cite{Cangiani2017Hp-VersionMeshes}.
Recent literature features several noteworthy implementations of polygonal discretizations in ALE schemes \cite{lipnikov2D,lipnikov2020ConservativeMeshes,Mazzia2020VirtualEnvironments,Gaburro2020HighChanges,Gaburro2021AChange,Wells2024AMethod}. These studies highlight the feasibility of achieving numerically stable solutions in ALE scenarios using various complex polygonal and polyhedral meshes. Additionally, some of these studies indicate the robustness of polygonal discretizations against time-dependent topological changes within ALE frameworks, even when faced with degenerate element faces \cite{Gaburro2020HighChanges,Mazzia2020VirtualEnvironments,Gaburro2021AChange,Wells2024AMethod}.

The Virtual Element Method (VEM) is a relatively new approach to discretising partial differential equations (PDEs) on polygonal and polytopic meshes \cite{basicprinciples,Ahmad2013EquivalentMethods,DaVeiga2014}. Extending the concepts of the Finite Element Method (FEM), the VEM expands the capabilities of the FEM by accommodating general polygonal or polyhedral elements. This adaptability makes it exceptionally suitable for scenarios involving complex geometries or non-standard mesh configurations. Initial applications of VEM have successfully addressed various types of PDEs, including elliptic, parabolic, and hyperbolic equations \cite{ellipticVEM,ConNonConVEM,Vacca2015VirtualMeshes,Vacca2017VirtualMeshes}. 
Moreover, the VEM has seen significant advancements in CFD, with applications ranging from Darcy flow to Navier-Stokes equations \cite{Antonietti2014AMeshes,daVeiga2017DivergenceMeshes,daVeiga2018VirtualMeshes,Liu2019TheEquations,Wang2019AMeshes,Irisarri2019StabilizedEquations,Adak2021AFormulation,Zhao2020TheProblem,Antonietti2022VirtualEquation,Verma2023VirtualMeshes}.

The use of the VEM for ALE schemes was first proposed in \cite{lipnikov2D} in which the VEM, specifically the polynomial projection operators, were employed to derive a conservative ALE scheme for the transformation of discrete data between two meshes. A velocity-based moving mesh algorithm was recently extended to the lowest order VEM for non-linear diffusion problems, attaining the same orders of convergence as the original finite element method \cite{Wells2024AMethod}.

Motivated by these results, this paper presents an Arbitrary Lagrangian-Eulerian Virtual Element Method (ALE-VEM) scheme for convection-diffusion equations. To extend the accuracy of this scheme beyond second-order, the isoparametric VEM of \cite{IsoVEM} is directly applied in the formulation. The proposed method is shown numerically to achieve high orders of accuracy in the $L^2$ and $H^1$ norms using the DUNE-VEM module \cite{Dedner2022ASpaces} which is part of the Distributed and Unified Numerics Environment \cite{Dedner2012Dune-Fem:Computing,Bastian2021TheDevelopments}.

The formulation is restricted to problems in which the ALE mapping is prescribed \textit{a priori} or can be solved for independently of the numerical solution. To demonstrate the extensibility of this method, this paper concludes with an extension of the velocity-based moving mesh VEM of \cite{Wells2024AMethod} to higher-order discretisations using the ALE-VEM framework. For this problem, the moving boundary is not prescribed \textit{a priori} and must be approximated at each discrete time level. Numerical results demonstrate that higher orders of convergence are achieved.

The layout of this paper is as follows. The convection-diffusion equation and a continuous ALE formulation are presented in Section \ref{sec::convection-diffusion}. The VEM framework is outlined in Section \ref{sec::VEM}. The isoparametric VEM and the discretisation of the ALE formulation are presented in Section \ref{sec::ALE-VEM}. Numerical convergence results are provided in Section \ref{sec::numerics}. In Section \ref{sec::VelocityBased}, the formulation and numerical results of the ALE-VEM scheme applied to the velocity-based moving mesh VEM are presented.

For an open bounded subset $\omega \subset \mathbb{R}^d$ where $d \ge 1$, the Sobolev space $W_p^m(\omega)$, for non-negative integer $m$ and $1\le p \le \infty$, is introduced with its norm $\norm{\cdot}_{m,p,\omega}$ and semi-norm $\abs{\cdot}_{m,p,\omega}$. In the special case of $p=2$, the Hilbert space $H^{m}(\omega)$ is defined with its norm $\norm{\cdot}_{m,\omega}$ and semi-norm $\abs{\cdot}_{m,\omega}$. $\mathbb{P}_k(\omega)$ denotes the space of polynomials of degree $k$ over $\omega$. 

The ALE-VEM is constructed on a bounded polygonal reference domain $\RefDomain \subset \rr^2$, with a corresponding reference coordinate system $\RefCoordinate = (\xi_1, \xi_2)$.
Functions or operators within the reference coordinate system are indicated using hat notation (e.g., $\hat{f}$, $\RefGrad \hat{f}$). Temporal dependency for functions, operators and bilinear forms are indicated through the use of a subscript $t$ or in the case of a discrete time level $t_n$ with a subscript of $n$.

\section{The Convection-diffusion Problem}\label{sec::convection-diffusion}

\subsection{The ALE Mapping}
We consider a bounded time-dependent domain $\Omega_t \subset \rr^2$ where $t \in [0,T]$ with a finite final time $T > 0$. 
The time-dependent ALE mapping is defined by $\ALEMap_t : \RefDomain \times [0,T] \rightarrow \Omega_t$,
and we assume for all time $t \in [0,T]$ that $\ALEMap_t$ is bi-Lipschitz, in line with assumptions made within ALE analysis \cite{Bonito2013Time-DiscreteStability,Bonito2013Time-discreteAnalysis,Gastaldi2001AElements,Formaggia1999StabilityElements} and that $ \ALEMap_t \in \left[ W^{k+1}_\infty (\RefDomain) \right]^2$ where $k \in \mathbb{N}$ denotes the degree of VEM being used. 
The Jacobian of the ALE mapping is given by
\[
    \Jacobian_t := \frac{\partial \ALEMap_t}{\partial \RefCoordinate} = \RefGrad \ALEMap_t,
\]
The determinant of $\Jacobian_t$  is uniformly bounded and strictly positive \cite{Bonito2013Time-DiscreteStability,Bonito2013Time-discreteAnalysis}; namely,  there exists a $\alpha \in (0,1)$ such that
\[
 \alpha \leq j_t := \det{\Jacobian_t}  \in L^{\infty}(\RefDomain) \qquad \forall t \in [0,T].
\]

\subsection{The Convection-diffusion Equation}
The linear convection-diffusion problem is defined as: for $t \in (0,
T]$, find $\rho(\PhysCoordinate,t)$ such that
\begin{align}
    \frac{\partial \rho}{\partial t} - \mu \Delta \rho + \nabla \cdot ( \mathbf{b} \rho ) &= f \qquad && \PhysCoordinate \in \Omega_t,\ t \in (0,T],\label{eq::A-DPDE} \\
    \rho &= \rho_0(\PhysCoordinate) &&\PhysCoordinate \in \Omega_0,\ t = 0, \nonumber \\
    \rho &= 0 &&\PhysCoordinate \in \partial \Omega_t,\ t \in [0,T], \nonumber
\end{align}
where $\b$ is a convective velocity field and $\mu \geq 0$ is a constant diffusivity parameter.

\subsection{Weak Formulation}
We define the space of functions which remain constant along domain trajectories by
\begin{equation}
    \mathcal{X}(\Omega_t) = \left\{ v:\ \Omega_t \rightarrow \rr\ :\ v = \hat{v} \circ \ALEMap^{-1}_t, \hat{v} \in H_0^1(\RefDomain) \right\}. \label{eq::TransportSpace}
\end{equation}
Any function $v \in \mathcal{X}(\Omega_t)$ has a vanishing material derivative with respect to the ALE velocity field \cite{Donea2004}
\begin{equation*}
    \frac{\partial v}{\partial t} + \ALEVelocity \cdot \nabla v = 0.
\end{equation*}
The Reynolds Transport Theorem \cite{Wesseling2001PrinciplesDynamics,Donea2004} states for $\rho$ and a given test function $v \in \mathcal{X}(\Omega_t)$ that
\begin{equation*}
    \frac{d}{dt} \int_{\Omega_t} \rho v\ d\PhysCoordinate = \int_{\Omega_t} \frac{\partial \rho v}{\partial t} + \nabla \cdot (\rho v \ALEVelocity)\ d\PhysCoordinate.
\end{equation*}
Application of the product rule and noting that the material derivative of $v$ vanishes provides
\begin{align*}
    \frac{d}{dt} \int_{\Omega_t} \rho v\ d\PhysCoordinate &= \int_{\Omega_t} v \left\{ \frac{\partial \rho }{\partial t} + \nabla \cdot (\rho \ALEVelocity) \right\} \ + \rho \left\{ \frac{\partial v }{\partial t} + \nabla \cdot \ALEVelocity \right\} \ d\PhysCoordinate \\
     &= \int_{\Omega_t} v \left\{ \frac{\partial \rho }{\partial t} + \nabla \cdot (\rho \ALEVelocity) \right\} \ d\PhysCoordinate.
\end{align*}
Substitution of the PDE \eqref{eq::A-DPDE} and integration by parts leads to
\begin{align*}
    \frac{d}{dt} \int_{\Omega_t} \rho v\ d\PhysCoordinate = \int_{\Omega_t} fv -\mu \nabla \rho \cdot \nabla v - (\ALEVelocity - \b) \cdot \rho \nabla v\  d\PhysCoordinate.
\end{align*}
The above equation can be written as follows
\begin{align}
    \frac{d}{dt} M_t(\rho,v) + \mu A_t(\rho,v) + B_t(\rho,v;\ALEVelocity) = l_t(v), \label{eq::AD-RTT-Form}
\end{align}
where,
\begin{align}
    M_t(\rho,v) &= \int_{\Omega_t} \rho v\ d\PhysCoordinate, \label{eq::M_tDef}\\
    A_t(\rho,v) &= \int_{\Omega_t} \nabla \rho \cdot \nabla v\ d\PhysCoordinate, \label{eq::A_tDef}\\
    B_t(\rho,v; \ALEVelocity) &= \int_{\Omega_t} (\ALEVelocity-\b) \cdot \rho \nabla v\ d\PhysCoordinate, \label{eq::B_tDef} \\
    l_t(v) &= \int_{\Omega_t} fv\ d\PhysCoordinate. \label{eq::l_tDef}
\end{align}

The stability of this formulation is well known under the assumption that $\nabla \cdot \b \in L^{\infty}(\Omega_t)$ and $f \in H^{-1}(\Omega_t)$ for all $t\in[0,T]$ \cite{Formaggia1999StabilityElements}.
The stability of the continuous formulation is independent on the choice of ALE mapping. Stability estimates have been presented for a FEM discretisation of the ALE method in which the stability is dependent on the choice of ALE mapping \cite{Formaggia1999StabilityElements,Formaggia2004StabilityALEFEM, Bonito2013Time-DiscreteStability}. In this paper we do not present any analysis on the ALE-VEM scheme, instead we assume that the choices of $\b$ and $f$ satisfy the aforementioned regularity conditions and that this is sufficient to propose a numerically stable method.

The continuous formulation given by Equation \eqref{eq::AD-RTT-Form} is independent of the choice of a conservative or non-conservative ALE approach in the continuous framework \cite{Formaggia1999StabilityElements}. As with FEMs, this is not the case for VEM discretisation of Equation \eqref{eq::AD-RTT-Form} and the VEM formulation in Section \ref{sec::ALE-VEM} is only applicable to a conservative ALE formulation of the convection-diffusion problem.

\section{The Virtual Element Method} \label{sec::VEM}
In this section, we outline the fundamental components of the VEM. The computational implementation of the VEM is well-documented and we omit the details in this paper \cite{DaVeiga2014,ConNonConVEM,Sutton2017TheMATLAB,Dedner2022ASpaces}.

\subsection{A Computational Reference Mesh}
We construct the VEM on a polygonal discretisation of $\RefDomain$ using the standard \textit{enhanced} VEM spaces \cite{Ahmad2013EquivalentMethods,DaVeiga2014,ConNonConVEM}. The polygonal mesh $\RefMesh$ of $\RefDomain$ is a collection of simple, non-overlapping polygonal reference elements $\RefE \in \RefMesh$.

The properties of a given $\RefE \in \RefMesh$ are defined as the diameter $h_{\RefE}$, barycentric coordinate $\PhysCoordinate_c^{\RefE} = (x_c^{\RefE}, y_c^{\RefE})$ and area $\abs{\RefE}$, boundary edges of $\RefE$ are denoted by $\hat{e} \subset \partial \RefE$. The global mesh size is defined by
\begin{equation*}
    h := \max_{\RefE \in \RefMesh} h_{\RefE}.
\end{equation*}

To ensure that optimal approximation results can be obtained, the following mesh regularity assumption is required \cite{basicprinciples}.
\begin{assumption}[Mesh Regularity] \label{assumption::MeshReg}~\\
Every $\RefE \in \RefMesh$ is a star-shaped domain or a finite union of star shaped domains with respect to a ball of radius greater than $\gamma h_{\RefE}$ for some uniform $\gamma > 0$. Additionally, for all edges $\hat{e} \in \partial \RefE$, the length of $\hat{e}$ is greater than $\delta h_{\RefE}$ for some  uniform $\delta > 0$.
\end{assumption}

\subsection{Polynomial Projection Operators}
The accuracy of the VEM is provided by using polynomial projection operators defined below \cite{basicprinciples,Ahmad2013EquivalentMethods,ellipticVEM}. 
For each of these projection operators there exists stability and accuracy results under the condition that Assumption \ref{assumption::MeshReg} is satisfied \cite{Brenner2008TheMethods,HHOBook}. 
\begin{definition}[The $\Pidel$ Operator] \label{def::Pidel}~\\
    The operator $\Pidel_k : H^1(\omega) \rightarrow \pp_{k}(\omega)$ is defined for all $v \in H^1(\omega)$ by
    \begin{align*}
        \int_{\omega} \nabla \Pidel_k v \cdot \nabla p\ d\PhysCoordinate &= \int_{\omega} \nabla v \cdot \nabla p\ d\PhysCoordinate \qquad \forall p \in \pp_k(\omega) \\
        \int_{\omega} v - \Pidel_k v\ d\PhysCoordinate &= 0.
    \end{align*}
\end{definition}

\begin{definition}[The $\PiO$ Operator] \label{def::Pi0}~\\
    The operator $\PiO_k: L^2(\omega) \rightarrow \pp_k(\omega)$ is defined for all $v \in L^2(\omega)$ by
    \begin{align*}
        \int_{\omega} \PiO_k v\ p\ d\PhysCoordinate &= \int_{\omega}  v\ p\ d\PhysCoordinate \qquad \forall p \in \pp_k(\omega).
    \end{align*}
\end{definition}

\begin{definition}[The $\PiI$ Operator] \label{def::Pi1}~\\
    The operator $\PiI_{k}: H^1(\omega) \rightarrow \pp_{k}(\omega)$ is defined for all $v \in H^1(\omega)$ by $\PiI_{k} v := \PiO_k \nabla v$ or equivalently
    \begin{align*}
        \int_{\omega} \PiI_{k} v \cdot \mathbf{p}\ d\PhysCoordinate &= \int_{\omega}  \nabla v \cdot \mathbf{p}\ d\PhysCoordinate \qquad \forall \mathbf{p} \in \left[ \pp_k(\omega) \right]^2.
    \end{align*}
\end{definition}

\subsection{Local \& Global VEM spaces}
The original VEM space is built on a given polygon $\RefE \in \RefMesh$ by considering a local boundary value problem subject to piecewise polynomial boundary conditions \cite{basicprinciples}. The boundary space is defined as
\begin{equation*}\label{eq::VEMBoundarySpace}
    \mathbb{B}_k(\partial \RefE) = \{ \hat{v}_h \in C^0(\partial \RefE)\ :\ \hat{v}_h|_e \in \mathbb{P}_k(\hat{e})\ \ \ \forall \hat{e} \subset \partial \RefE \},
\end{equation*}
from which the original local VEM space of degree $k$ can then be defined as 
\begin{equation*}\label{eq::VEMSpaceOriginal}
    W_k(\RefE) = \{ \hat{v}_h \in H^1(\RefE)\ :\ \hat{v}_h|_{\partial \RefE} \in \mathbb{B}_k(\partial \RefE),\ \Delta \hat{v}_h|_{E} \in \pp_{k-2}(\RefE) \}.
\end{equation*}

In order to develop the VEM for problems beyond Poisson's Equation, the local VEM space had to be modified such that the full $L^2$ projection was computable. A solution was presented in \cite{Ahmad2013EquivalentMethods} and is commonly used as the default VEM space in most recent discretisations. 

The approach to enhance the VEM space is to use the gradient projection operator $\Pi^\nabla$ to supplement the remaining internal moments required to compute the $L^2$ projection. The original local VEM space is enlarged via,
\begin{equation*}
    \tilde{W}_k(\RefE) = \{ \hat{v}_h \in H^1(\RefE)\ :\  \hat{v}_h|_{\partial \RefE} \in \mathbb{B}_k(\partial \RefE),\ \Delta \hat{v}_h \in \mathbb{P}_k(\RefE) \},
\end{equation*}
and then restricted to define the enhanced VEM space,
\begin{equation}\label{eq::ClassicalLocalVEMSpace}
    V_k(\RefE) = \left\{ \hat{v}_h \in \tilde{W}_k(\RefE)\ :\ \int_{\RefE} (\hat{v}_h - \Pi^{\nabla}_k \hat{v}_h) \hat{q}\ d\PhysCoordinate = 0\ \ \ \forall \hat{q} \in \mathbb{P}_k(\RefE) \backslash \pp_{k-2}(\RefE) \right\}.
\end{equation}

The global VEM space of degree $k$ is defined as
\begin{equation}\label{eq::classicGlobalVEMSPace}
    \hat{V}_h = \left\{ \hat{v}_h \in H^1(\RefDomain) \:\ \hat{v}_h|_{\hat{E}} \in V_k(\RefE)\ \forall \RefE \in \RefMesh \right\}.
\end{equation}
Using the VEM to solve PDEs with homogeneous Dirichlet boundary conditions requires a restricted global VEM space with zero trace. This is denoted as $\hat{V}_{h,0} = \hat{V}_h \cap H_0^1(\RefDomain)$.

\subsection{Degrees of Freedom}
The construction of a VEM relies on choosing a set of degrees of freedom (DoFs). These DoFs serve two purposes in the method: they provide the necessary information to construct computable polynomial approximation operators, and they uniquely identify a virtual element function on any given polygon, which would otherwise be an unknown solution to a local boundary value problem.

The space of scaled monomials of degree $k$ on a given element $\RefE \in \RefMesh$ is defined by
\begin{equation*}\label{eq::ScaledMonomialSpace}
\mathcal{M}_k(\RefE) = \left\{ \left( \frac{x - x_c}{h_{\RefE}} \right)^{\alpha} \left( \frac{y-y_c}{h_{\RefE}} \right)^{\beta}\ :\ \alpha,\beta \in \mathbb{Z}{\geq 0},\ \alpha+\beta = k \right\},
\end{equation*}
from which the VEM DoFs are defined as follows.
\begin{definition}[Degrees of Freedom for the Virtual Element Method]\label{def::VEMDoFs}~\\
Let $\RefE\in \RefMesh$ and $V_k(\RefE)$ be the local VEM space defined in Equation \eqref{eq::ClassicalLocalVEMSpace}. The degrees of freedom of a given function $\hat{v}_h \in V_k(\RefE)$ are defined as follows:
\begin{itemize}
\item The point values of $\hat{v}_h$ at each vertex of $\RefE$.
\item The point values of $\hat{v}_h$ at the Gauss-Lobatto quadrature points of order $k-1$ on each edge $\hat{e} \subset \partial \RefE$.
\item The internal moments of $\frac{1}{|\RefE|}\int_{\RefE} \hat{v}_h \hat{q}\ dx$ for all $\hat{q}\in \mathcal{M}_{k-2}(\RefE)$.
\end{itemize}
\end{definition}
A proof that these constitute a unisolvent set of DoFs is provided in \cite{basicprinciples}. The projections $\Pi^\nabla_k$ and $\Pi^0_k$ are computable using only these degrees of freedom \cite{Ahmad2013EquivalentMethods,DaVeiga2014}.

\subsection{Representation of VEM Functions}
For a given VEM function $\hat{v}_h \in \hat{V}_h$, the DoFs are denoted by $\text{dof}_i(\hat{v}_h)$ for $i=1,\ldots,\Ndofs$, with $\Ndofs$ being the number of DoFs on $\hat{V}_h$. The Lagrangian VEM basis function is introduced as
\begin{align}\label{eq::VEMBasisFun}
\left\{ \hat{\varphi}_i \right\}_{i=1}^{\Ndofs} \subset \hat{V}_h, \qquad \text{dof}_i(\hat{\varphi}_j) = \delta_{i,j} \text{ for } i,j=1,\ldots,\Ndofs.
\end{align}
VEM functions and interpolants can be defined respectively using the DoFs and basis functions as
\begin{align} \label{eq::VEMLagrangeExpansion}
\hat{v}_h = \sum_{i=1}^{\Ndofs} \text{dof}_i(\hat{v}_h) \hat{\varphi}_i,\ \ \   \hat{v}_I = \sum_{i=1}^{\Ndofs} \text{dof}_i(\hat{v}) \hat{\varphi}_i \qquad \forall \hat{v}_h \in \hat{V}_h,~ \hat{v} \in  L^1(\RefDomain) \cap C^0(\RefDomain).
\end{align}

\subsection{Approximation of Linear and Bilinear Forms}
Let $\mathcal{A}(\cdot,\cdot): H^1_0(\RefDomain) \times H^1_0(\RefDomain) \rightarrow \rr$ and $l(\cdot): H^1_0(\RefDomain) \rightarrow \rr$ be a symmetric bilinear and linear form respectively. Consider the general variational problem: find $\hat{\rho} \in H^1_0(\RefDomain)$ such that $\mathcal{A}(\hat{\rho},\hat{v}) = l(\hat{v})$ for all $\hat{v} \in H^1_0(\RefDomain)$. 
The novelty of the VEM is the introduction of discrete approximations $\mathcal{A}_h(\cdot,\cdot)$ and $l_h(\cdot)$ that use the DoFs and the projection operators to approximate the integral equations to a sufficient degree of accuracy. These approximations are constructed via local element contributions
\begin{align*}
    \mathcal{A}_h(\hat{\rho}_h,\hat{v}_h) = \sum_{\RefE \in \RefMesh} \mathcal{A}_h^{\RefE}(\hat{\rho}_h,\hat{v}_h), \qquad l_h(\hat{v}_h) = \sum_{\RefE\in \RefMesh} l_h^{\RefE}(\hat{v}_h) \qquad \forall \hat{v}_h \in V_{h,0}.
\end{align*}
To enforce the coercivity of the discrete bilinear form, a stabilisation term is introduced $S^{\RefE}(\cdot,\cdot)$ to ensure that the kernel of $\mathcal{A}_h^{\RefE}(\cdot,\cdot)$ scales like the kernel of $\mathcal{A}^{\RefE}(\cdot,\cdot)$ \cite{basicprinciples,ConNonConVEM}. We provide further details on these discrete forms in Section \ref{sec::ALE-VEM}.
\section{A Conservative ALE-VEM Formulation}\label{sec::ALE-VEM}

\subsection{Approximating the ALE Mapping}
In this section, we present the conservative ALE discretisation where the ALE mapping and velocity are known \textit{a priori} and interpolated component-wise into the tensor product VEM space $\left[ \hat{V}_h \right]^2$ (c.f. Equation \eqref{eq::VEMLagrangeExpansion}):
\begin{align}
\ALEMap_{h,t} = \left[ \ALEMap_{I,t}^x,\ \ALEMap_{I,t}^y \right], \qquad
\hat{\ALEVelocity}_{h,t} = \left[ \hat{\ALEVelocity}_{I,t}^x,\ \hat{\ALEVelocity}_{I,t}^y \right],
\end{align}
where the $x,y$ superscripts denote the $x$ and $y$ components of $\ALEMap_t$ and $\hat{\ALEVelocity}_t$ respectively.
These functions, being defined on the reference domain, depend spatially on $\RefCoordinate$ and have a time-dependent set of degrees of freedom. The formulation of this ALE scheme does not require that $\ALEMap_{h,t}$ is the VEM interpolant of $\ALEMap_t$. Alternatively, the $\ALEMap_{h,t}$ only needs to be a sufficiently accurate approximation of $\ALEMap_t$ at a given time $t \in [0,T]$ for the isoparametric VEM \cite{IsoVEM}. In particular, the approach of the moving mesh method of Section \ref{sec::VelocityBased} approximates the ALE mapping at each discrete time-level.

\subsection{Moving Virtual Element Spaces}
The virtual domain generated by the ALE map $\ALEMap_{h,t}$ is denoted by $\Omega_{h,t}$ and the corresponding virtual mesh of transformed elements is denoted by $\mathcal{T}_{h,t}$.
The discrete counterpart to Equation \eqref{eq::TransportSpace} is defined using the discrete ALE mapping
\begin{align}
    \mathcal{X}_h (\Omega_{h,t}) = \left\{ v: \Omega_t \rightarrow \rr\ :\ v = \hat{v} \circ \ALEMap_{h,t}^{-1},\ \hat{v} \in \hat{V}_h \right\}. \label{eq::DiscreteTransportSpace}
\end{align}
For the VEM we choose a basis to be a subset of $\mathcal{X}_h$
\begin{align}
    \left\{ \varphi_i(\mathbf{x},t)\ :\ \varphi_i(\mathbf{x},t) = \hat{\varphi}_i(\RefCoordinate) \circ \ALEMap_{h,t}^{-1} \ \forall t \in [0,T] \right\}_{i=1}^{\Ndofs} \subset \mathcal{X}_h(\Omega_{h,t}), \label{eq::MovingVEMBasis}
\end{align}
with $\hat{\varphi}_i$ being the canonical VEM basis function of $\hat{V}_{h}$ defined using Equation \eqref{eq::VEMBasisFun}. Using this basis, we define a time-dependent VEM function by 
\begin{equation}
      v_{h,t}(\mathbf{x},t) = \sum_{i=1}^{\Ndofs} \text{dof}_i(v_h(\mathbf{x},t)) \varphi_i(\mathbf{x},t), \label{eq::MovingVEMFunction}
\end{equation}
where $\varphi_i$ is defined by Equation \eqref{eq::MovingVEMBasis}. We emphasis that the DoFs of the discrete function $v_{h,t}$ in Equation \eqref{eq::MovingVEMBasis} are time-dependent so in general $v_{h,t} \notin \mathcal{X}_h(\Omega_{h,t})$, instead $v_{h,t}$ is a time-dependent linear combination of elements of $\mathcal{X}_h(\Omega_{h,t})$. The moving VEM space can then be defined, using Equations \eqref{eq::MovingVEMBasis} and \eqref{eq::MovingVEMFunction}, as
\begin{equation}
    V_{h,t} = \left\{ v_{h,t} \in H^1(\Omega_{h,t})\ :\ v_{h,t} = \hat{v}_{h,t} \circ \ALEMap_{h,t}^{-1},\ \hat{v}_{h,t} \in \hat{V}_{h} \right\}. \label{eq::MovingVEMSPace}
\end{equation}
To impose the homogeneous Dirichlet boundary conditions, the restriction of these VEM spaces to zero boundary conditions are defined by $V_{h,t,0} := V_{h,t} \cap H^1_0(\Omega_{h,t})$ and $\mathcal{X}_{h,0} := \mathcal{X}_{h} \cap H^1_0(\Omega_{h,t})$ respectively.

\subsection{The Isoparametric VEM}
To achieve higher-order accuracy with the VEM solution we require a higher-order representation of the moving domain. An isoparametric VEM allows for the representation of the ALE map using the same space $\hat{V}_h$. In this work, we consider the first isoparametric VEM of \cite{IsoVEM} which transforms the variational problem from the discrete time-dependent domain, represented by a VEM function, onto the computational reference domain (see also \cite{lipnikov2D,lipnikov2020ConservativeMeshes,Bachini2021Arbitrary-orderSurfaces} for similar approaches to domain transformations). 

Firstly, we consider at a fixed time $t \in [0,T]$ approximating the semi-discrete ALE formulation of Equation \eqref{eq::AD-RTT-Form} on the VEM approximation of the moving domain $\Omega_{h,t} \approx \Omega_t$. An approximation of the Jacobian operator is introduced to approximate a change of variables from $\Omega_{h,t}$ to $\RefDomain$. For a given time $t \in [0,T]$ and reference element $\hat{E} \in \RefMesh$, we define this operator approximation and it's corresponding determinant by
\begin{align}
    \Jacobian_{h,t} = \Pi^1_{k-1} \ALEMap_{h,t}, \qquad j_{h,t} = \det{\Jacobian_{h,t}}. \label{eq::DiscreteJacobian}
\end{align}
Using these polynomial, we approximate the transformation of the gradient of a VEM function $v_h \in V_{h,t}$ from a physical element $E_h \in \mathcal{T}_{h,t}$ to the corresponding reference element $\hat{E} \in \RefMesh$ and the determinant of the Jacobian via
\begin{align}
    \nabla v_h \approx \Jacobian_{h,t}^{-\top} \Pi^1_{k-1} \RefGrad \hat{v}_h, \qquad \det{\Jacobian_t} \approx j_{h,t}.\label{eq::DiscreteGradient}
\end{align}
Analysis on the accuracy of these approximations is provided in \cite{IsoVEM}. In particular, it is shown that for sufficiently small mesh size $h$ that $j_{h,t} > 0$ for all $t \in [0,T]$ and we use this assumption in the formulations of this method.

Important to the formulation is the assumption that the PDE data, in this case $\mu$, $\mathbf{b}$ and $f$ are available and can be computed to a sufficient degree of accuracy on the computational reference domain. In the numerical experiments of Sections \ref{sec::numerics} and \ref{sec::VelocityBased}, the PDE data will either be a VEM function or explicitly known on $\RefDomain$.


\subsection{A Semi-discretisation}\label{sec::Semi}
In this VEM the numerical solution of the convection-diffusion equation is given by $\rho_{h,t} \in V_{h,t,0}$. This solution is only implicitly known in the ALE coordinates as the ALE-VEM scheme is computed using only the representation of $\rho_{h,t}$ and the test functions $v_{h,t} \in \mathcal{X}_{h}(\Omega_{h,t})$ in the reference coordinates. 

For ease of reading, we drop the temporal subscript for the solution and test functions in the following formulations, instead the time-dependency will be described by the temporal subscript in the bilinear forms of the method. 

The ALE-VEM semi-discrete formulation is given as follows: for a given $t \in (0,T]$ find $\rho_{h,t} \in V_{h,t,0}$ such that 
\begin{equation*}
    \frac{d}{dt} M_{h,t}(\rho_h,v_h) + \mu A_{h,t}(\rho_h,v_h) + B_{h,t}(\rho_h,v_h;\ALEVelocity_h) = l_{h,t}(v_h) \qquad \forall v_h \in \mathcal{X}_{h,0}(\Omega_{h,t}).
\end{equation*}

The VEM discretisations of Equations \eqref{eq::M_tDef}, \eqref{eq::A_tDef}, \eqref{eq::B_tDef} and \eqref{eq::l_tDef} are given by element-wise contributions on the reference mesh
\begin{align*}
    M_{h,t}(\rho_h,v_h) = \sum_{\hat{E} \in \RefMesh} M_{h,t}^{\hat{E}}(\rho_h,v_h),\qquad A_{h,t}(\rho_h,v_h) = \sum_{\hat{E} \in \RefMesh} A_{h,t}^{\hat{E}}(\rho_h,v_h), \\
    B_{h,t}(\rho_h,v_h;\ALEVelocity_h) = \sum_{\hat{E} \in \RefMesh} B_{h,t}^{\hat{E}}(\rho_h,v_h;\ALEVelocity_h), \qquad l_{h,t}(v_h) = \sum_{\hat{E} \in \RefMesh} l_{h,t}^{\hat{E}}(v_h).
\end{align*}
The local element contributions are defined using Equations \eqref{eq::DiscreteJacobian} and \eqref{eq::DiscreteGradient} by
\begin{align}
    M_{h,t}^{\hat{E}}(\rho_h,v_h) &= \int_{\hat{E}} \Pi^0_k \hat{\rho}_h\ \Pi^0_k \hat{v}_h\ j_{h,t}\ d\RefCoordinate + h_{\hat{E}}^2 S^{\hat{E}}(\hat{\rho}_h - \Pi^0_k \hat{\rho}_k, \hat{v}_h - \Pi^0_k \hat{v}_k), \label{eq::M_ht} \\
    A_{h,t}^{\hat{E}}(\rho_h,v_h) &= \int_{\hat{E}} \Jacobian_{h,t}^{-\top} \Pi^1_{k-1} \hat{\rho}_h\ \Jacobian_{h,t}^{-\top} \Pi^1_{k-1} \hat{v}_h\ j_{h,t}\ d\RefCoordinate \nonumber \\
    &\hspace{4cm} + S^{\hat{E}}(\hat{\rho}_h - \Pi^0_k \hat{\rho}_k, \hat{v}_h - \Pi^0_k \hat{v}_k), \label{eq::A_ht} \\
    B_{h,t}^{\hat{E}}(\rho_h,v_h; \ALEVelocity_h) &= \int_{\hat{E}} (\Pi^0_k \hat{\ALEVelocity}_h- \hat{\b}) \cdot \Pi^0_k \rho_h \Jacobian_{h,t}^{-\top} \Pi^1_{k-1} \hat{v}_h\ j_{h,t}\ d\RefCoordinate, \label{eq::B_ht} \\
    l_{h,t}^{\hat{E}}(v_h) &= \int_{\hat{E}} \hat{f}\ \Pi^0_k \hat{v}_h\ j_{h,t}\ d\RefCoordinate, \label{eq::l_ht}
\end{align}
where $S^{\hat{E}}(\cdot,\cdot)$ is the standard dofi-dofi stabilisation term \cite{basicprinciples,Ahmad2013EquivalentMethods,DaVeiga2014}
\begin{equation}
    S^{\RefE}(\hat{z}_h,\hat{v}_h) = \sum_{i=1}^{\Ndofs} \text{dof}_i(\hat{z}_h) \cdot \text{dof}_i(\hat{v}_h) \qquad \forall \hat{z}_h, \hat{v}_h \in \hat{V}_{h,t}.
\end{equation}

The analysis of \cite{IsoVEM} implies that each bilinear form approximation is $O(h^k)$ accurate for a sufficiently small reference mesh size $h$, for any $t \in [0,T]$, provided the reference mesh is shape-regular and the assumptions on the ALE mapping detailed in Section \ref{sec::convection-diffusion} hold.

\subsection{A fully Discrete Scheme}\label{sec::thetaScheme}
A method of lines approach is taken to perform the integration of the weak formulation with respect to time.
We define the fully discrete ALE-VEM scheme using the $\theta$-scheme. For $\theta \in [0,1]$ and $0 \leq t_n < t_{n+1} \leq t_{N^t} = T$ with $\Delta t_{n+1} = t_{n+1}-t_n$, the time derivative of $M_{h,t}(\rho_h,v_h)$ is approximated via
\begin{align}
    \frac{M_{h,n+1}(\rho_{h},v_h) - M_{h,n}(\rho_h,v_h)}{\Delta t} &= \theta l_{h,n+1}(\rho_h,v_h) + (1-\theta)l_{h,n}(\rho_h,v_h) \nonumber \\
    &-\mu \left(  \theta A_{h,n+1}(\rho_h,v_h) + (1-\theta) A_{h,n}(\rho_h,v_h) \right) \nonumber\\
    &- \left( \theta B_{h,n+1}(\rho_h,v_h; \ALEVelocity_h) + (1-\theta) B_{h,n}(\rho_h,v_h; \ALEVelocity_h) \right). \label{eq::MassThetaScheme}
\end{align}
By rearranging the terms of Equation \eqref{eq::MassThetaScheme}, we define the matrices $\boldsymbol{\mathcal{G}}_{n+1}$ and $\boldsymbol{\mathcal{H}}_{n}$ by
\begin{align*}
(\mathcal{G}_{n+1})_{i,j} &= M_{h,n+1}(\varphi_i,\varphi_j) + \Delta t \theta [ A_{h,n+1}(\varphi_i,\varphi_j) + B_{h,n+1}(\varphi_i,\varphi_j; \ALEVelocity_{h}) ] \\
(\mathcal{H}_n)_{i,j} &= M_{h,n}(\varphi_i,\varphi_j) + \Delta t (\theta - 1) [A_{h,n}(\varphi_i,\varphi_j) + B_{h,n}(\varphi_i,\varphi_j;\ALEVelocity_{h})] ,
\end{align*}
and the vector $\boldsymbol{\mathcal{F}}_n$ is defined as
\begin{align*}
    (\mathcal{F}_n)_i &= \Delta t \theta l_{h,n+1}(\varphi_i) + \Delta t (1-\theta) l_{h,n}(\varphi_i).
\end{align*}

The fully discrete ALE-VEM scheme is the defined as follows: given the DoFs vectors of $\boldsymbol{\rho}_n$, $\ALEVelocity_{n+1}$ and $\ALEVelocity_{n}$, find $\boldsymbol{\rho}_{n+1}$ such that
\begin{align}
        \boldsymbol{\mathcal{G}}_{n+1} \boldsymbol{\rho}_{n+1} = \boldsymbol{\mathcal{H}}_n \boldsymbol{\rho}_n + \boldsymbol{\mathcal{F}}_n. \label{eq::ALE-VEMFullyDisc}
\end{align}
\section{Numerical Experiments}\label{sec::numerics}

\subsection{Computational Reference Mesh}
The Centroidal Voronoi Tessellation (CVT) mesh structure is used to test the implementation of the ALE-VEM scheme \cite{Senechal1995SpatialDiagrams,Du2006ConvergenceTessellations}. A sequence of randomly generated CVT meshes is used for each test in this section. Each mesh in the sequence is sampled such that the mesh size $h$ roughly halves with each refinement. The mesh files are generated using PolyMesher \cite{Talischi2012PolyMesher:Matlab} within MATLAB and imported into DUNE. The reference domain is taken as the unit square $\RefDomain = [0,1]^2$.

\subsection{Error Computation}
We assess the numerical error in the solution at the final time $t_N = T$. The Solution $H^1$ and $L^2$ errors are approximated by the discrete norms of $\norm{\cdot}_{h,1}$ and $\norm{\cdot}_{h,0}$ respectively. These are defined as
\begin{align*}
    \norm{\rho_{h,N} - \rho_{N}}_{h,1}^2 &:= \sum_{\hat{E} \in \RefMesh} \int_{\hat{E}} \abs{ \Jacobian_{h,N}^{-\top} \ \Pi^1_{k-1} ( \hat{\rho}_{h,N} - \hat{\rho}_N ) }^2\ j_{h,N}\ d\RefCoordinate, \\
    \norm{\rho_{h,N} - \rho_{N}}_{h,0}^2 &:= \sum_{\hat{E} \in \RefMesh} \int_{\hat{E}} \abs{ \ \Pi^0_{k} ( \hat{\rho}_{h,N} - \hat{\rho}_N ) }^2\ j_{h,N}\ d\RefCoordinate.
\end{align*}
This choice of discrete norm is an approximation of the error on the physical domain. Discrete error norms on the reference domain can also be considered, in line with the error estimate of \cite{IsoVEM}. Computing errors on the reference domain produced the same empirical orders of convergence (EOC) presented in this section and are not presented here.

\subsection{Experiment Conditions}
The numerical experiments are performed within the DUNE software environment \cite{Bastian2021TheDevelopments, Dedner2012Dune-Fem:Computing, Alkamper2016TheModule} using the Python bindings presented in \cite{Dedner2020PythonModule}. Unless stated otherwise, the Crank-Nicolson method ($\theta=0.5$) is used to integrate the fully discrete system \eqref{eq::ALE-VEMFullyDisc} between discrete time-levels. 
The time step size $\Delta t$ is reduced in all convection-diffusion simulations according to $\Delta t^2 \sim h^{k+1}$ such that the expected orders of convergence in the spatial discretisation are produced. All numerical experiments in this section are run to a final time of $T=0.01$.

To ensure that the accuracy of the ALE-VEM scheme can be assessed, problems with inhomogenous boundary conditions have to be considered. These conditions are enforced by applying a Dirichlet boundary condition on the method at each time step using the interpolant of the DoFs of the true solution on the moving boundary
\begin{equation*}
    \hat{\rho}_{h,n} |_{\partial \RefDomain} = \hat{\rho}_{I,n} \ \ \ \ n=0,...,N.
\end{equation*}
The DUNE-UFL library is used to symbolically compute the forcing data and the explicit expressions of these terms are omitted in this paper. 

\subsection{A Time-independent ALE Map}
In the first test we validate the numerical method by considering a time-independent domain transformation. The CE and Warped Square mappings of \cite{lipnikov2D}. These are respectively defined as
\begin{align}
    \ALEMap(\RefCoordinate) &= \left[ \xi_1 +  \xi_1 \xi_2 (1-\xi_1)/2,\ \xi_2 + \xi_1 \xi_2 (1-\xi_2)/2 \right], \label{eq::CEMap} \\
    \ALEMap(\RefCoordinate) &= \left[ \sin{\frac{\xi_1 \pi}{3}},\ e^{\xi_2} \right]. \label{eq::WS}
\end{align}
We choose the diffusivity parameter to be $\mu=1$ and the convective velocity to be $\mathbf{b} = \mathbf{x}$. The forcing data $f$ and Dirichlet boundary conditions are chosen such that the true solution of the convection-diffusion problem in the physical domain is given as
\begin{align*}
    \rho(x,y,t) = \exp(-\pi^2 t) \sin(\pi x) \sin(\pi y).
\end{align*}
Numerical results for both mappings are given in Table \ref{table::time-independent::CE} and Table \ref{table::time-independent::WS}. Here we observe the expected convergence rates in the $H^1$ and $L^2$ norms of $O(h^k)$ and $O(h^{k+1})$ respectively when run to a final time of $T=0.01$.

\begin{table}[ht]
\centering
\caption{$L^2$ and $H^1$ solution error data for the convection-diffusion equation on the time-independent CE ALE mapping \eqref{eq::CEMap} for $k=1,2,3$.}\label{table::time-independent::CE}
\begin{tabular}{|c|c|c|c|c|c|c|}
\hline
\( h \) & \multicolumn{2}{c|}{\( k=1 \)} & \multicolumn{2}{c|}{\( k=2 \)} & \multicolumn{2}{c|}{\( k=3 \)} \\
\cline{2-7}
& L2  & H1 & L2  & H1 & L2  & H1 \\
\hline
0.1041 & 2.604e-3 & 3.349e-2 & 3.798e-4 & 8.826e-3 & 8.583e-5 & 2.302e-3 \\
0.0480 & 8.283e-4 & 2.219e-2 & 6.426e-5 & 2.184e-3 & 7.870e-6 & 3.276e-4 \\
0.0231 & 2.549e-4 & 1.090e-2 & 7.560e-6 & 4.257e-4 & 6.000e-7 & 4.754e-5 \\
0.0113 & 5.968e-5 & 5.222e-3 & 7.800e-7 & 9.109e-5 & 4.000e-8 & 5.580e-6 \\
\hline
rate & 2.038 & 1.073 & 3.104 & 2.027 & 3.873 & 2.889 \\
\hline
\end{tabular}
\end{table}

\begin{table}[ht]
\centering
\caption{$L^2$ and $H^1$ solution error data for the convection-diffusion equation on the time-independent Warped Square ALE mapping \eqref{eq::WS} for $k=1,2,3$.}\label{table::time-independent::WS}
\begin{tabular}{|c|c|c|c|c|c|c|}
\hline
\( h \) & \multicolumn{2}{c|}{\( k=1 \)} & \multicolumn{2}{c|}{\( k=2 \)} & \multicolumn{2}{c|}{\( k=3 \)} \\
\cline{2-7}
& L2  & H1 & L2  & H1 & L2  & H1 \\
\hline
0.1041 & 1.981e-2 & 1.718e-1 & 2.808e-3 & 3.527e-2 & 2.169e-4 & 4.519e-3 \\

0.0480 & 6.246e-3 & 7.194e-2 & 3.680e-4 & 6.361e-3 & 1.176e-5 & 6.429e-4 \\

0.0231 & 1.444e-3 & 2.343e-2 & 5.507e-5 & 1.146e-3 & 7.100e-7 & 8.500e-5 \\
0.0113 & 3.708e-4 & 9.798e-3 & 6.970e-6 & 2.274e-4 & 5.000e-8 & 1.083e-5 \\
\hline
rate & 1.964 & 1.362 & 2.986 & 2.330 & 3.841 & 2.930 \\
\hline
\end{tabular}
\end{table}

\subsection{A Pure Diffusion Problem}
For this experiment, a heat equation problem is considered in which there is a zero convection term $\mathbf{b} = 0$, leading to a linear diffusion problem on a moving domain. We set $\mu =1$ and $f = 0$ and choose an initial conditions such that the solution is given by
\begin{align*}
    \rho(x,y,t) = \exp(-2 \pi^2 t) \sin(\pi x) \sin(\pi y).
\end{align*}
The domain is transformed by the linearised CE mapping from Equation \eqref{eq::CEMap} 
\begin{align}
    \ALEMap_t(\RefCoordinate) = \RefCoordinate +  \frac{t}{2T} \xi_1 \xi_2  \left[ 1 - \xi_1, 1 - \xi_2 \right],\label{eq::CELin}
\end{align}
where $T=0.01$ is the final time in the simulation. The vorticial motion (VM) map \cite{lipnikov2D} is defined as the solution of the system of ODEs:
\begin{align}
    \dot{x} &= 2 \sin (\pi x) \cos(\pi y), \label{eq::VM::dx} \\
    \dot{y} &= -2 \cos(\pi x)\sin (\pi y). \label{eq::VM::dy}
\end{align}
 In the case of the VM mapping defined by Equations \eqref{eq::VM::dx} and \eqref{eq::VM::dy}, only the velocity-field of the ALE mapping is provided and the non-linear nature of the velocity field requires an explicit time integration scheme. To match the second-order accuracy of the Crank-Nicolson time-stepping scheme, we employ Heun's method (modified Euler \cite{Burden2016NumericalEd.}) to integrate Equations \eqref{eq::VM::dx} and \eqref{eq::VM::dy} over time. 
The numerical results are presented in Tables \ref{table::CE} and \ref{table::VM} for the CE and VM ALE mappings respectively. Again we observe $O(h^{k+1})$ and $O(h^k)$ orders of convergence in the $L^2$ and $H^1$ norms for both mappings.

\begin{table}[ht]
\centering
\caption{$L^2$ and $H^1$ solution error data for the convection-diffusion equation on the CE ALE mapping \eqref{eq::CELin} for $k=1,2,3$.}\label{table::CE}
\begin{tabular}{|c|c|c|c|c|c|c|}
\hline
\( h \) & \multicolumn{2}{c|}{\( k=1 \)} & \multicolumn{2}{c|}{\( k=2 \)} & \multicolumn{2}{c|}{\( k=3 \)} \\
\cline{2-7}
& L2  & H1 & L2  & H1 & L2  & H1 \\
\hline
0.1041 & 9.922e-3 & 1.038e-1 & 3.980e-3 & 5.173e-2 & 3.545e-3 & 4.611e-2 \\
0.0480 & 2.144e-3 & 3.212e-2 & 5.433e-4 & 8.260e-3 & 2.194e-4 & 2.887e-3 \\
0.0231 & 5.182e-4 & 1.294e-2 & 6.833e-5 & 1.497e-3 & 1.371e-5 & 2.159e-4 \\
0.0113 & 1.230e-4 & 6.017e-3 & 7.650e-6 & 2.794e-4 & 8.600e-7 & 2.087e-5 \\
\hline
rate  & 2.085    & 1.115    & 3.157    & 2.413    & 3.986    & 3.366    \\
\hline
\end{tabular}
\end{table}

\begin{table}[ht]
\centering
\caption{$L^2$ and $H^1$ solution error data for the convection-diffusion equation on the VM ALE mapping \eqref{eq::VM::dx}, \eqref{eq::VM::dy} for $k=1,2,3$.}\label{table::VM}
\begin{tabular}{|c|c|c|c|c|c|c|}
\hline
\( h \) & \multicolumn{2}{c|}{\( k=1 \)} & \multicolumn{2}{c|}{\( k=2 \)} & \multicolumn{2}{c|}{\( k=3 \)} \\
\cline{2-7}
& L2  & H1 & L2  & H1 & L2  & H1 \\
\hline
0.1041 & 1.142e-2 & 1.179e-1 & 7.431e-4 & 8.660e-3 & 6.926e-5 & 9.192e-4 \\
0.0480 & 2.778e-3 & 3.891e-2 & 1.184e-4 & 1.378e-3 & 4.300e-6 & 1.145e-4 \\
0.0231 & 7.015e-4 & 1.042e-2 & 1.555e-5 & 2.283e-4 & 2.700e-7 & 1.455e-5 \\
0.0113 & 1.709e-4 & 3.462e-3 & 1.760e-6 & 4.918e-5 & 2.000e-8 & 1.710e-6 \\
\hline
rate  & 2.033    & 1.552    & 3.132    & 2.215    & 3.760    & 3.074    \\
\hline
\end{tabular}
\end{table}

\subsection{A General Convection-diffusion Problem}
\begin{figure}[ht]
\centering
\begin{minipage}[t]{.45\textwidth}
  \centering
  \includegraphics[width=1\linewidth]{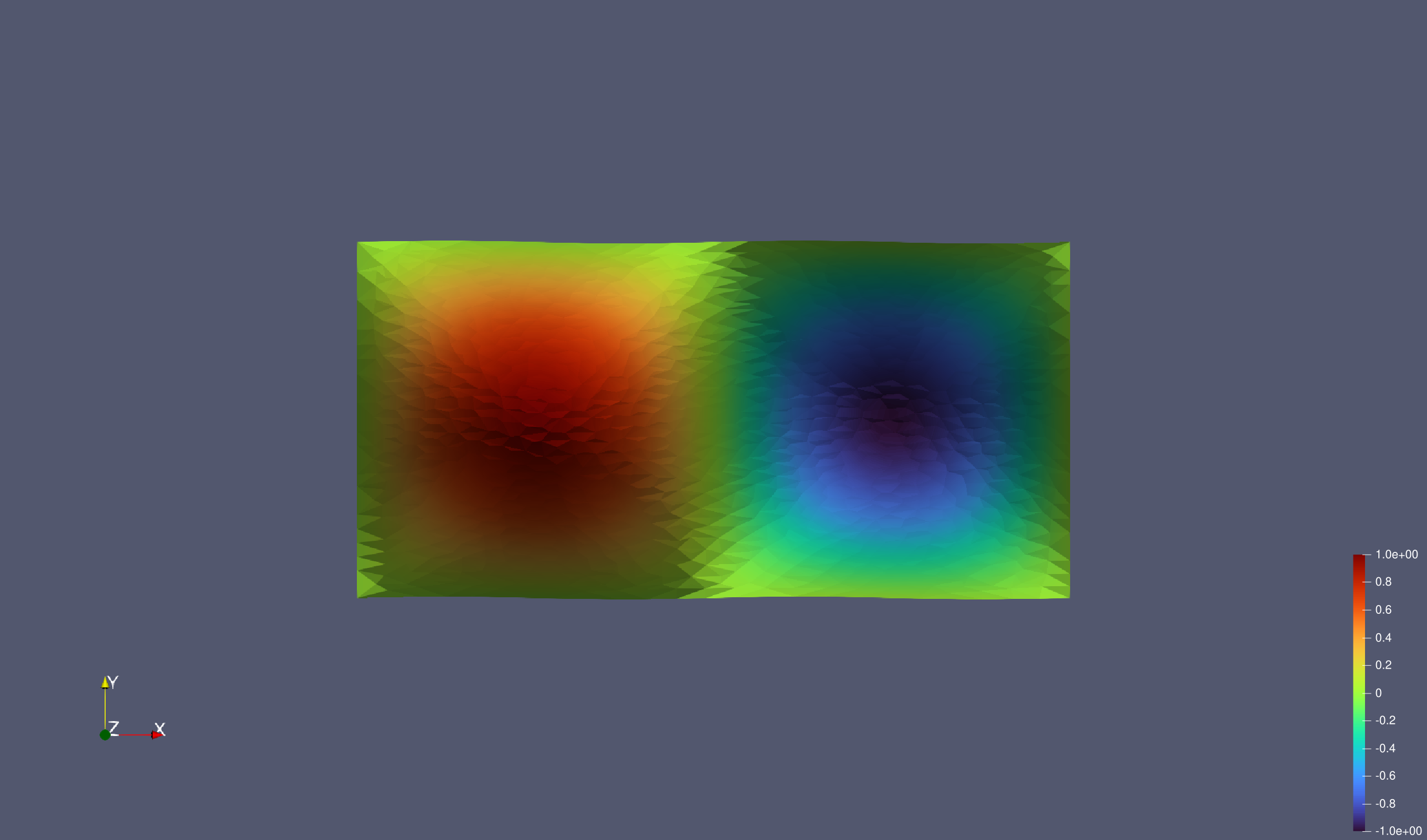}
\end{minipage}
~
\begin{minipage}[t]{.45\textwidth}
  \centering
  \includegraphics[width=1\linewidth]{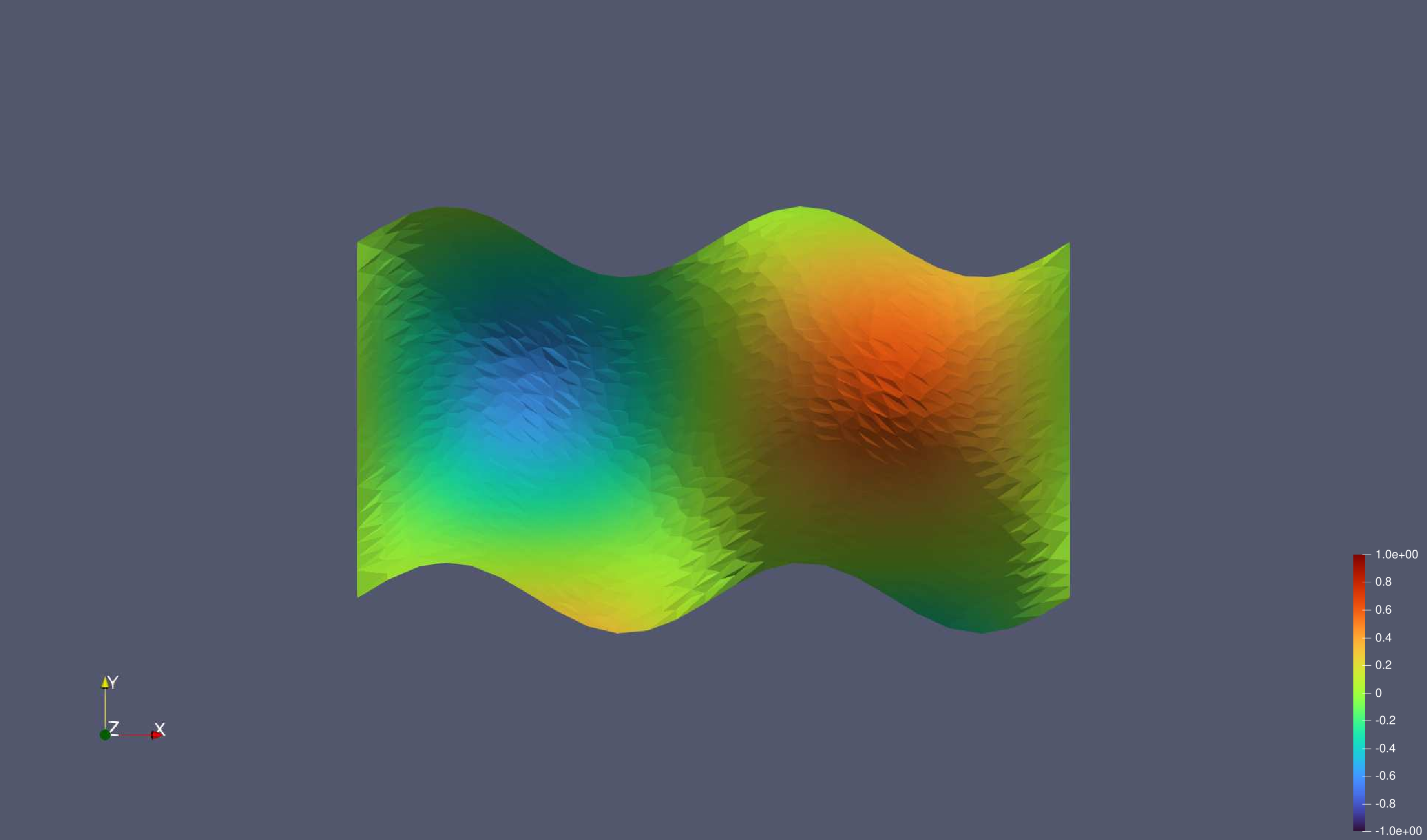}
\end{minipage}
\\
\begin{minipage}[t]{.45\textwidth}
  \centering
  \includegraphics[width=1\linewidth]{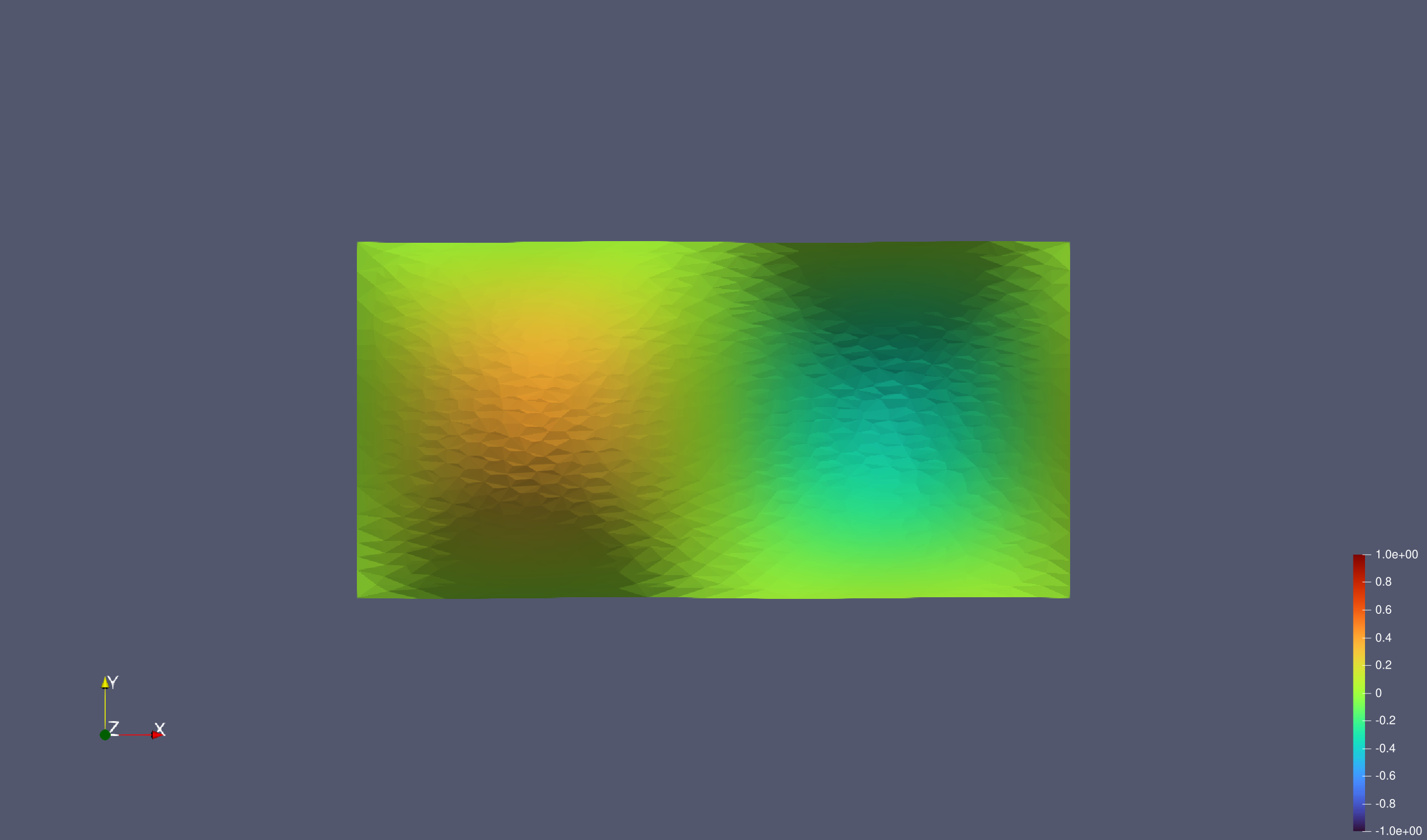}
\end{minipage}
~
\begin{minipage}[t]{.45\textwidth}
  \centering
  \includegraphics[width=1\linewidth]{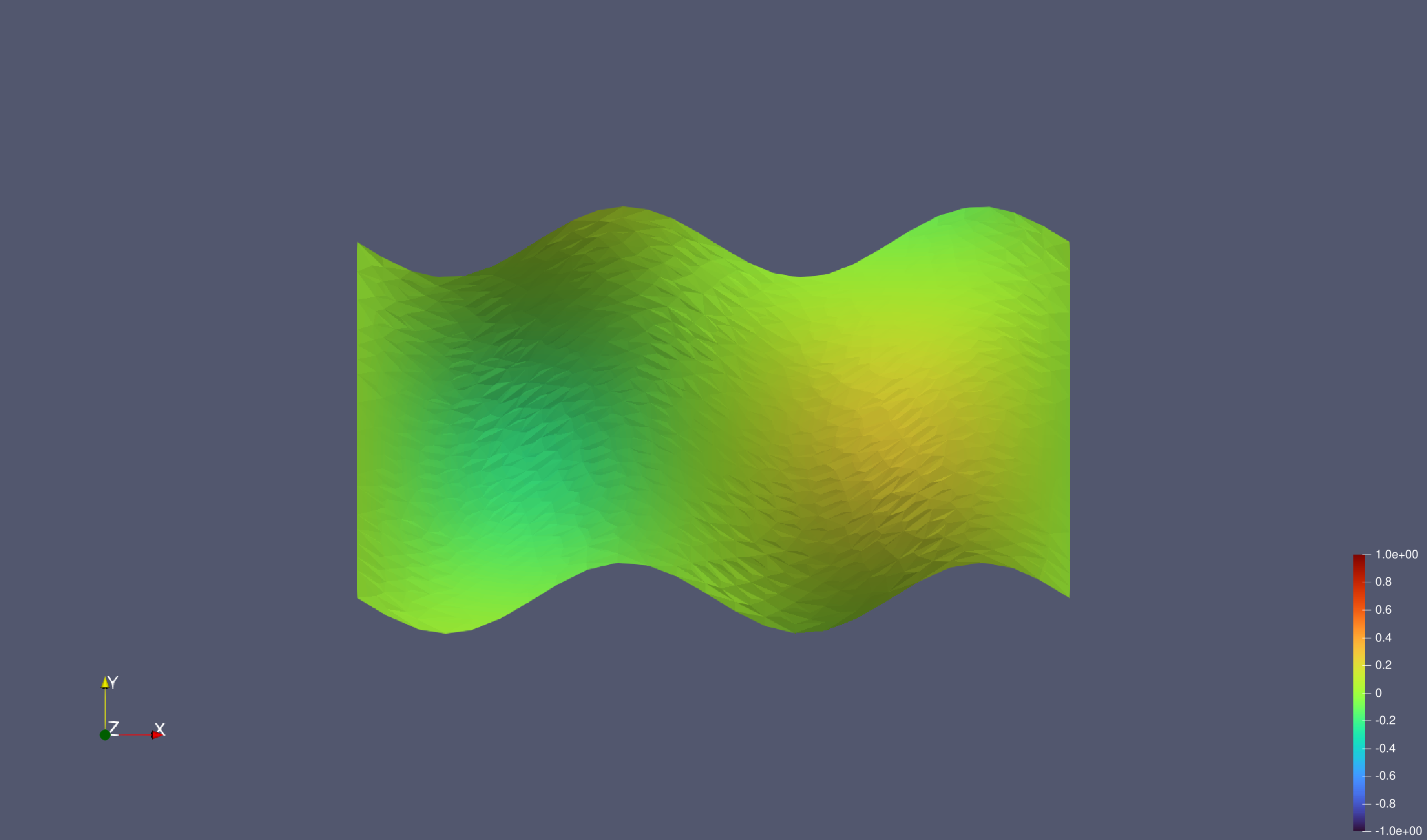}
\end{minipage}
\captionof{figure}{Solution snapshots of the convection-diffusion equation using a quadratic VEM and a reference mesh of 800 elements. The snapshots are taken at times $t=0$ (top left), $t=0.025$ (top right), $t=0.05$ (bottom left) and $t=0.075$ (bottom right).}
\label{fig::ADMesh}
\end{figure}

\begin{figure}[ht]
\centering
\begin{minipage}[t]{.45\textwidth}
  \centering
  \includegraphics[width=1\linewidth]{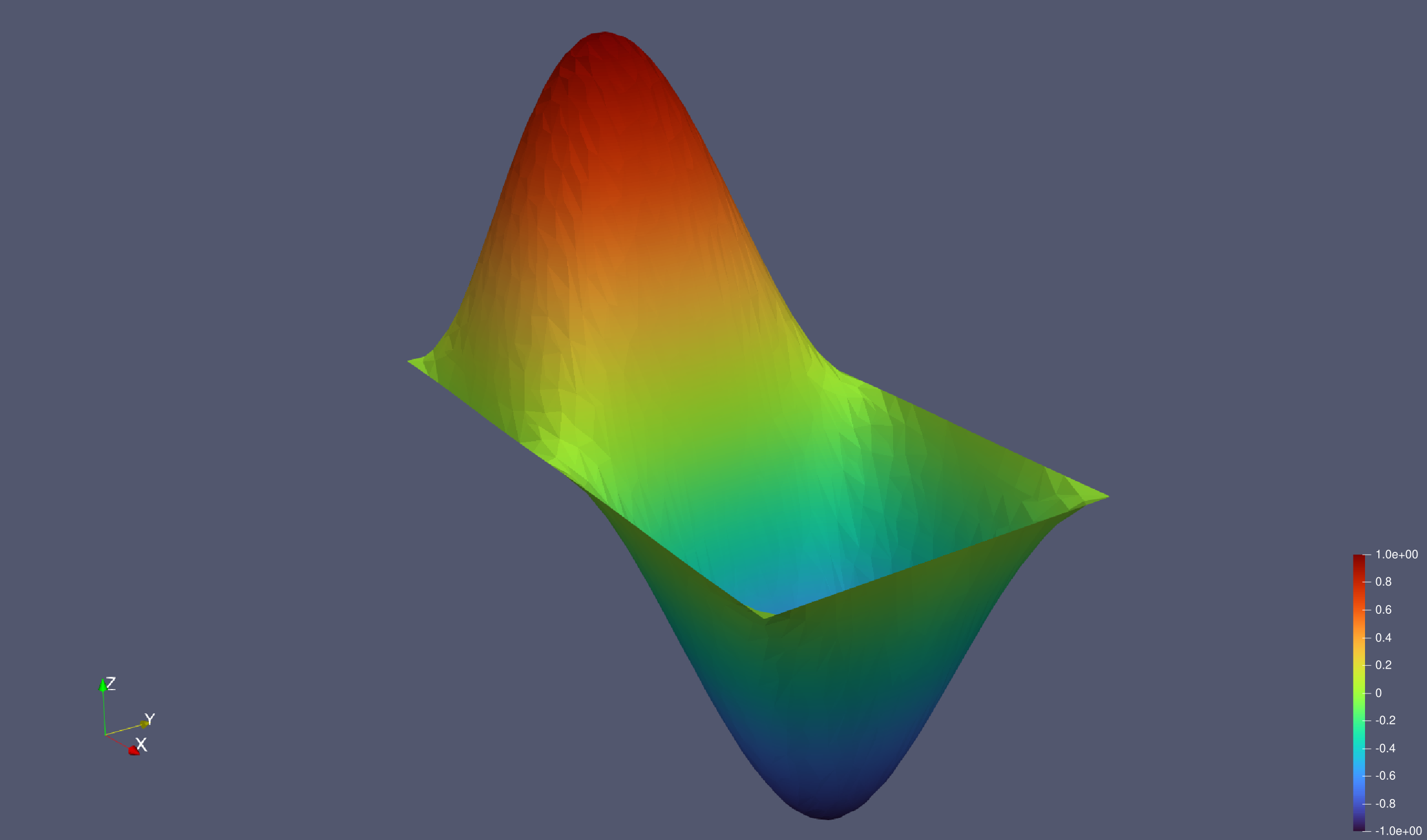}
\end{minipage}
~
\begin{minipage}[t]{.45\textwidth}
  \centering
  \includegraphics[width=1\linewidth]{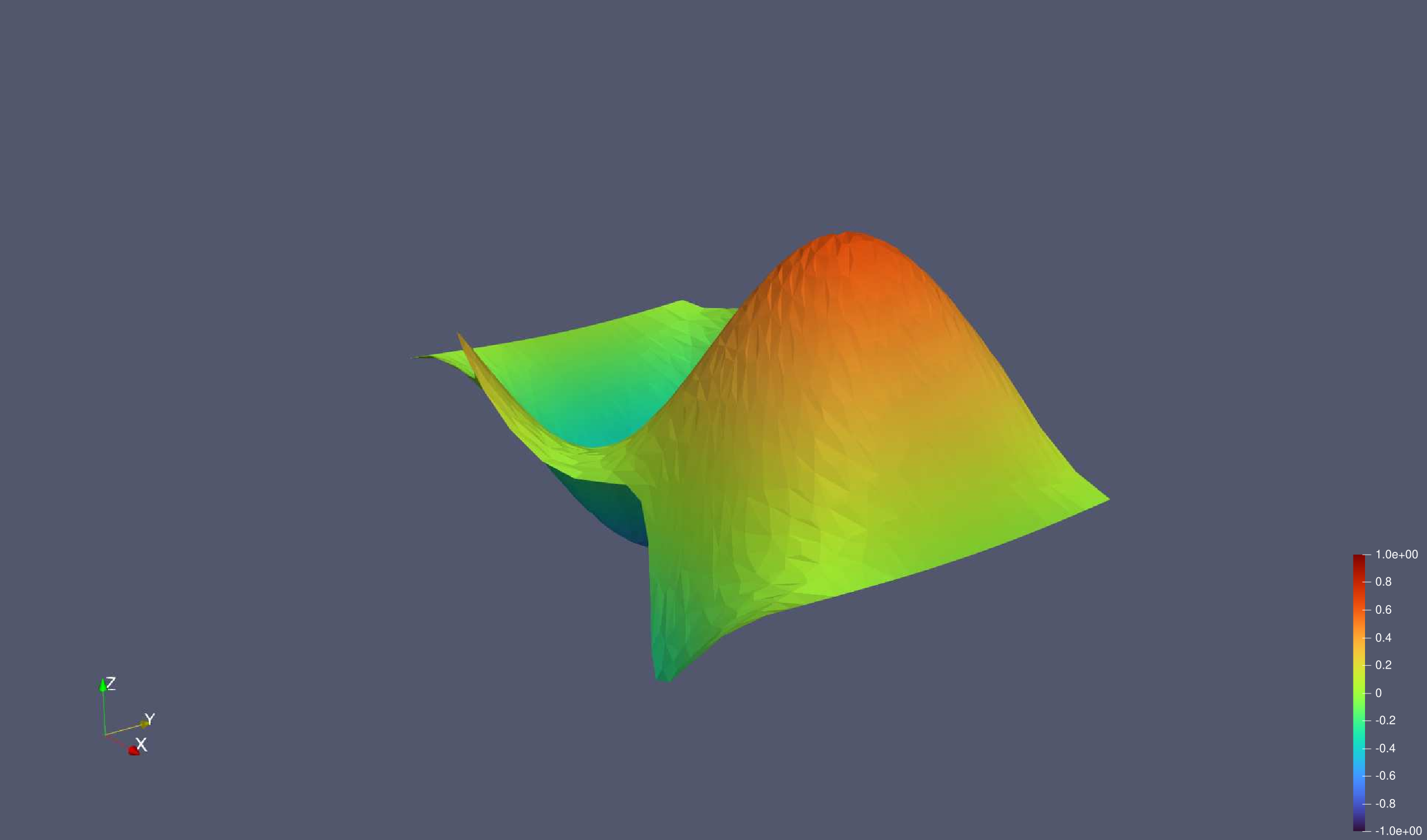}
\end{minipage}
\\
\begin{minipage}[t]{.45\textwidth}
  \centering
  \includegraphics[width=1\linewidth]{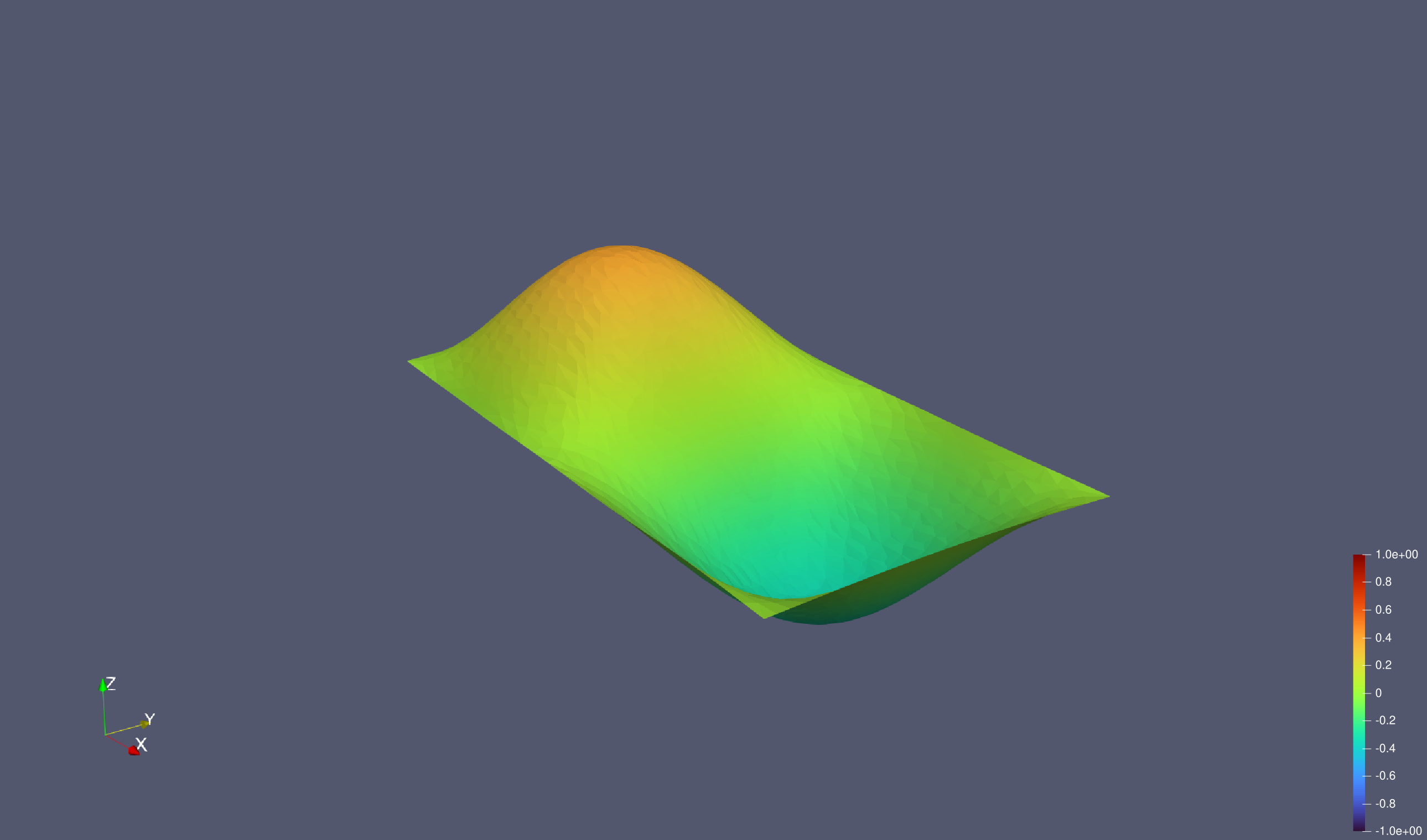}
\end{minipage}
~
\begin{minipage}[t]{.45\textwidth}
  \centering
  \includegraphics[width=1\linewidth]{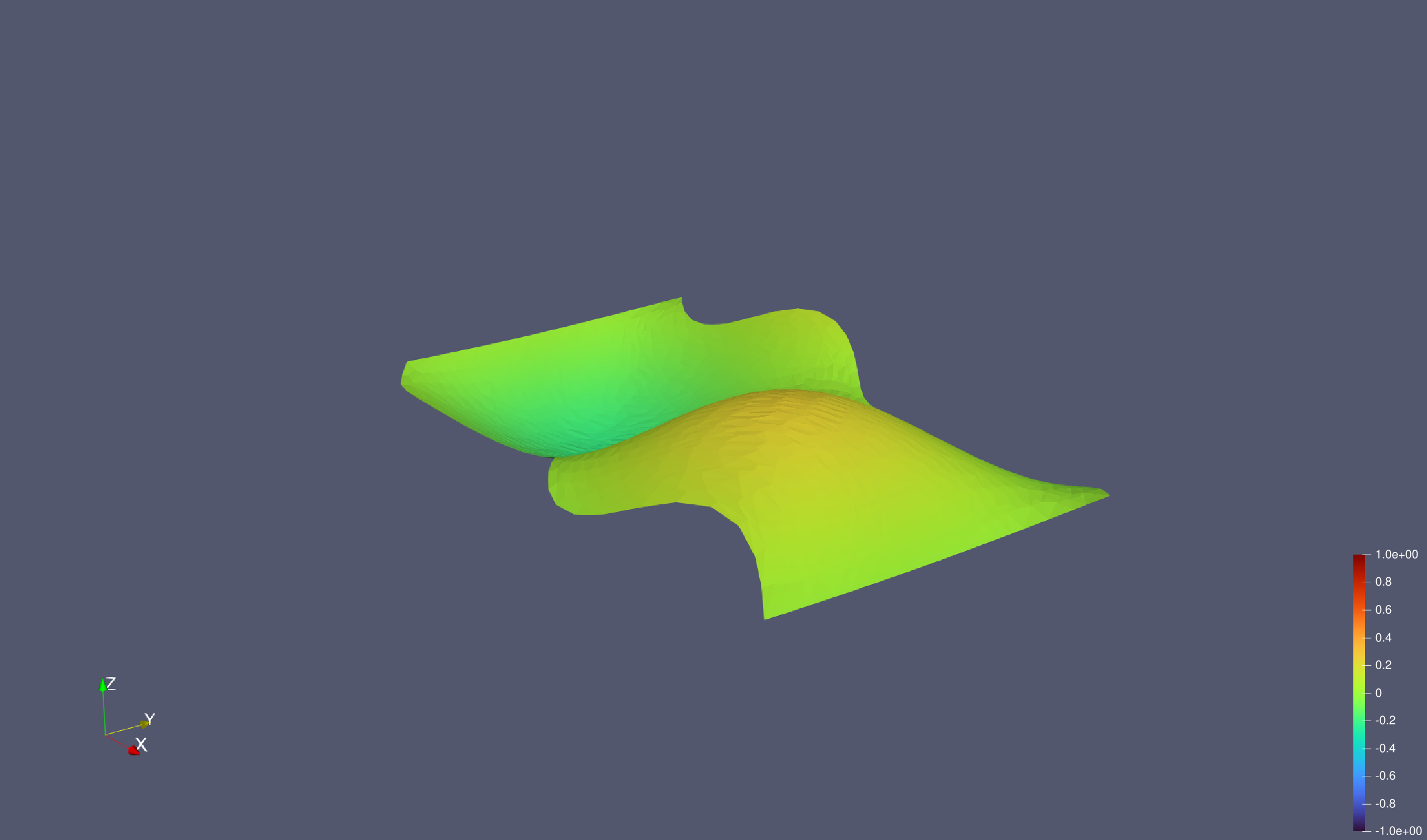}
\end{minipage} 
\captionof{figure}{Solution snapshots of the convection-diffusion equation using a quadratic VEM and a reference mesh of 800 elements. The snapshots are taken at times $t=0$ (top left), $t=0.025$ (top right), $t=0.05$ (bottom left) and $t=0.075$ (bottom right).}
\label{fig::ADSol}
\end{figure}

In this experiment, we consider the numerical example of a solution to the convection-diffusion equation with a travelling feature. 
The domain is transformed from $\RefDomain = [0,1]^2$ to a rectangular domain $[0,2] \times [0,1]$ which oscillates in the y direction. This transformation is defined as
\begin{align}
    \ALEMap_{t}(\RefCoordinate) = \RefCoordinate + [ \xi_1,\ A \sin(\pi u_y t) ],\label{eq::OSMap}
\end{align}
where $A >0$ and $u_y$ are user specified parameters which control the amplitude and frequency of the oscillations respectively. The ALE velocity field of this transformation is given as
\begin{align}
    \ALEVelocity = [0, A u_y \pi \cos(\pi u_y t)].
\end{align}
The parameters of the convection-diffusion Equation \eqref{eq::A-DPDE} are set as $\mu =1$, $\mathbf{b} = \mathbf{x}$ and $f$ is chosen such that the true solution of the PDE is 
\begin{align*}
    \rho(\mathbf{x},t) = \exp(-\pi^2 t) \sin (\pi (y - Asin(\pi u_y t)))\sin (\pi (x -u_x t)),
\end{align*}
where $u_x$ is another user specified parameter that controls the speed of travel of the solution in the $x$ direction. In our experiment we consider the ALE mapping given by $A=1/10$ and $u_x = u_y = 20$ and run the simulation to a final time of $T=0.01$. From Table \ref{table::OS}, we observe the expected orders of convergence in the $L^2$ and $H^1$ norms.
Moving mesh and solution profile snapshots are provided in Figures \ref{fig::ADMesh} and \ref{fig::ADSol} respectively. 

\begin{table}[ht]
\centering
\caption{$L^2$ and $H^1$ solution error data for the convection-diffusion equation on the ALE mapping \eqref{eq::OSMap} for $k=1,2,3$.}\label{table::OS}
\begin{tabular}{|c|c|c|c|c|c|c|}
\hline
\( h \) & \multicolumn{2}{c|}{\( k=1 \)} & \multicolumn{2}{c|}{\( k=2 \)} & \multicolumn{2}{c|}{\( k=3 \)} \\
\cline{2-7}
& L2  & H1 & L2  & H1 & L2  & H1 \\
\hline
0.1041 & 3.049e-2 & 2.722e-1 & 8.848e-3 & 1.072e-1 & 6.000e-3 & 5.045e-2 \\
0.0480 & 7.946e-3 & 1.008e-1 & 9.124e-4 & 2.521e-2 & 3.744e-4 & 6.153e-3 \\
0.0231 & 2.094e-3 & 3.930e-2 & 9.607e-5 & 5.194e-3 & 2.339e-5 & 8.432e-4 \\
0.0113 & 5.146e-4 & 1.702e-2 & 1.075e-5 & 1.063e-3 & 1.460e-6 & 1.080e-4 \\
\hline
rate  & 2.028    & 1.278    & 3.077    & 2.289    & 4.003    & 2.937    \\
\hline
\end{tabular}
\end{table}

\section{A Velocity-based Moving Mesh Method}\label{sec::VelocityBased}
To conclude this paper, we extend the conservative ALE-VEM to a velocity-based moving mesh method \cite{Baines2005,Baines2006,MMFEMBoundaries,Baines2009AProblems,Marlow2011MovingEquations,Baines2011}. The lowest-order VEM was successfully applied and tested to non-linear diffusion problems \cite{Wells2024AMethod}. In this work, we extend the simplest test case of a similarity solution to the Porous Medium Equation (PME) using quadratic and cubic virtual elements.

\subsection{The Porous Medium Equation}
The PME is a second-order non-linear diffusion problem that admits a family of similarity solutions with compact support. Modelling the support of these solutions as a moving boundary problem makes the PME an ideal test case for a moving mesh method. Extensive examples and theoretical analysis of the PME is provided in \cite{Vazquez2006TheEquation} and numerical properties of this solution are discussed in \cite{Baines2011,Ngo2017AEquation}. We consider a simple quasi-linear case of the PME, which in relation to Equation \eqref{eq::A-DPDE}, has  a non-linear convection term \( \mathbf{b} = - \nabla \rho \) and a zero diffusion parameter and forcing term \( \mu = f = 0\). In addition to the homogeneous boundary condition, we require an additional zero-flux boundary condition across the moving boundary \cite{Baines2005}. The PME is given as: for $t \in (0,T]$, find $\rho(\mathbf{x},t)$ such that 
\begin{align*}
    \frac{\partial \rho}{\partial t} &= \nabla \cdot (\rho \nabla \rho ) \qquad && \PhysCoordinate \in \Omega_t, t \in (0,T], \\
    \rho &= \rho_0(\PhysCoordinate) \qquad && \PhysCoordinate \in \Omega_0, \\
    \rho &= 0 \qquad && \PhysCoordinate \in \partial \Omega_t, t \in [0,T], \\
    \rho \nabla \rho \cdot \mathbf{n}_t &= 0 \qquad && \PhysCoordinate \in \partial \Omega_t, t \in [0,T],
\end{align*}
where $\mathbf{n}_t$ denotes the time-dependent outward unit normal vector of the moving boundary $\partial \Omega_t$.

\subsection{Solving for the ALE Velocity Field}
We derive a weak formulation of the velocity field on the following two assumptions of mass conservation and an irrotational ALE velocity field for all $t \in [0,T]$
\begin{align*}
    \frac{d}{dt} M_{t}(\rho,v) = 0,\ \forall v \in \mathcal{X}(\Omega_t), \qquad \exists \phi \in H^1(\Omega_t)  \text{ s.t. } \ALEVelocity = \nabla \phi.
\end{align*}
Under these assumptions we can derive a weak formulation at a fixed point in time for the velocity potential \cite{Baines2011}: given $\rho \in H^1(\Omega_t)$, find $\phi \in H^1(\Omega_t)$ such that $a(\phi, v) = d(v)$ for all $v \in H^1(\Omega_t)$, where 
\begin{align*}
    a(\phi,v) = \int_{\Omega_t} \rho \nabla \phi \cdot \nabla v\ d\PhysCoordinate, \qquad d(v) = -\int_{\Omega_t} \rho \nabla \rho \cdot \nabla v\ d\PhysCoordinate. 
\end{align*}
The ALE velocity field is then reconstructed by an $L^2$ projection of the gradient of $\phi$: given $\phi \in H^1(\Omega_t)$, find $\ALEVelocity \in \left[ H^1(\Omega_t) \right]^2$ such that $m(\ALEVelocity,v) = b(v)$ for all $v \in H^1(\Omega_t)$ where
\begin{align*}
    m(\ALEVelocity,v) = \int_{\Omega_t} \ALEVelocity\ v\ d\PhysCoordinate, \qquad b(v) = \int_{\Omega_t} \nabla \phi\ v\ d\PhysCoordinate.
\end{align*}

We discretise the potential and velocity reconstruction equations using the isoparametric VEM \cite{IsoVEM}.
Given $\rho_h \in V_{h,n,0}$ at time $t_n \in [0,T]$, find $\phi_h \in V_{h,n}$ such that $a_{h,n}(\phi_h,v_h) = d_{h,n}(v_h)$ for all $v_h \in \mathcal{X}_{h,n}$ where
\begin{align*}
    a_{h,n}(\phi_h,v_h) &= \sum_{\hat{E} \in \RefMesh} \int_{\hat{E}} \Pi^0_k \hat{\rho}_h \ \Jacobian^{-\top}_{h,n} \Pi^1_{k-1} \hat{\phi}_h \cdot \Jacobian^{-\top}_{h,n} \Pi^1_{k-1} \hat{v}_h\ j_{h,n}\ d\RefCoordinate \\
    &\hspace{4cm} + \bar{\rho}_h S^{\hat{E}} \left( \hat{\phi}_h - \Pi^0_k \hat{\phi}_h, \hat{v}_h - \Pi^0_k \hat{v}_h \right), \\
    d_{h,n}(v_h) &= - \sum_{\hat{E} \in \RefMesh} \int_{\hat{E}} \Pi^0_k \hat{\rho}_h\ \Jacobian_{h,n}^{-\top} \Pi^1_{k-1} \hat{\rho}_h \cdot \Jacobian_{h,n}^{-\top} \Pi^1_{k-1} \hat{v}_h\ j_{h,n}\ d\RefCoordinate.
\end{align*}
Given $\phi_h \in V_{h,n}$ at time $t_n \in [0,T]$, find $\ALEVelocity_h \in \left[ V_{h,n} \right]^2$ such that $m_{h,n}(\ALEVelocity_h,\mathbf{v}_h) = b_{h,n}(\mathbf{v}_h)$ for all $\mathbf{v}_h \in \left[ \mathcal{X}_{h,n} \right]^2$ where
\begin{align*}
    m_{h,n}(\ALEVelocity_h,\mathbf{v}_h) &= \sum_{\hat{E} \in \RefMesh}\int_{\hat{E}} \Pi^0_{k} \hat{\ALEVelocity}_h \cdot \Pi^0_{k} \hat{\mathbf{v}}_h\ j_{h,n}\ d\RefCoordinate \\
    &\hspace{4cm}+h_{\hat{E}}^2 S^{\hat{E}} \left( \hat{\ALEVelocity}_h - \Pi^0_k \hat{\ALEVelocity}_h, \hat{\mathbf{v}}_h - \Pi^0_k \hat{\mathbf{v}}_h \right), \\
    b_{h,n}(\mathbf{v}_h) &= \sum_{\hat{E} \in \RefMesh} \int_{\hat{E}} \Jacobian_{h,n}^{-\top} \Pi^1_{k-1} \hat{\phi}_h \cdot \Pi^0_k \hat{\mathbf{v}}_h\ j_{h,n}\ d\RefCoordinate. 
\end{align*}

\subsection{The ALE Equation}
In the interest of brevity, we only consider in this paper the case where $\ALEVelocity_h$ approximates only the Lagrangian velocity field of the domain. Examples using alternative ALE velocity fields have been studied in other works \cite{Marlow2011MovingEquations,Baines2011,Wells2024AMethod}. The Forward Euler time-stepping scheme ($\theta = 0$) is used for both the moving mesh and the ALE update, in line with existing implementations of this moving mesh method. Consequently, the ALE Equation \eqref{eq::MassThetaScheme} simplifies for the PME to
\begin{align}
    \frac{M_{h,n+1}(\rho_{h},v_h) - M_{h,n}(\rho_h,v_h)}{\Delta t} &= - B_{h,n}(\rho_h,v_h; \ALEVelocity_h), \label{eq::MassThetaSchemePME}
\end{align}
Noting that the definition of $B_{h,n}(\cdot,\cdot;\cdot)$ differs slightly from Equation \eqref{eq::B_ht} due to the non-linear convection term of the PME
\begin{equation}
    B_{h,n}(\rho_h,v_h;\ALEVelocity_{h}) = \sum_{\hat{E} \in \RefMesh} \int_{\hat{E}} \Pi^0_k \hat{\rho}_h\ \Jacobian_{h,n}^{-\top} \Pi^1_{k-1} \hat{v}_h \cdot \left( \Jacobian_{h,n}^{-\top} \Pi^1_{k-1} \hat{\rho}_h + \Pi^0_k \hat{\ALEVelocity}_h \right)\ j_{h,t}\ d\RefCoordinate. \nonumber
\end{equation}
We can couple Equation \eqref{eq::MassThetaSchemePME} with a Forward Euler time-stepping scheme for computing the ALE map to complete the moving mesh method
\begin{align*}
    M_{h,n+1}(\rho_{h},v_h) &=  M_{h,n}(\rho_h,v_h) - \Delta t B_{h,n}(\rho_h,v_h; \ALEVelocity_h) \qquad &&\forall v_h \in V_{h,n,0},  \\
    \ALEMap_{h,n+1} &= \ALEMap_{h,n} + \Delta t \ALEVelocity_{h,n}.
\end{align*}
A more comprehensive summary of the structure of this moving mesh method is provided in \cite{Baines2011}.

\subsection{Numerical Experiments}
The velocity-based moving mesh method is tested against a similarity solution of the PME. There exists a family of radially symmetric solutions on a given initial circular domain of radius $r_0$ for the PME defined in \cite{Wesseling2001PrinciplesDynamics} and given by
\begin{equation}\label{eq::PMEsolution}
\rho(r,t) = 
\begin{cases}
\frac{1}{\lambda(t)^2} \left( 1 - \left( \frac{r}{r_0 \lambda(t)} \right)^2 \right) \ \ \ &|r| \leq r_0\lambda(t) \\
0\ \ \ &\text{otherwise}
\end{cases},
\end{equation}
where $d$ is the spatial dimension, $r_0$ is the initial radius, and
\begin{equation}\label{eq::PMEextras}
\lambda(t) = \left( \frac{t}{t_0} \right)^{\frac{1}{4}},\ \ \ t_0 = \frac{r_0^2}{8}. \nonumber
\end{equation}

In order to achieve higher than second order accuracy, we require a higher order representation of the initial domain $\Omega_{0}$ \cite{WellsThesis,IsoVEM}. For linear finite element and virtual element approaches, a polygonal interpolant of the initial domain is sufficient to attain $O(h^2)$ accuracy in the solution of the PME \cite{Baines2011,Wells2024AMethod}. The computational domain is defined as the polygonal interpolant of $\Omega_0$ and we define $\mathbf{g}: \partial \RefDomain \rightarrow \partial \Omega_0$ as a Lipschitz mapping between the interpolated and true initial boundary. The initial ALE transformation is then approximated by solving a harmonic extension of $\mathbf{g}$ into $\RefDomain$ \cite{Donea2004,Richter2017}
\begin{align*}
    -\Delta \ALEMap_{0} &= \mathbf{0} \qquad && \RefCoordinate \in \RefDomain, \\
    \ALEMap_0 &= \mathbf{g} \qquad && \RefCoordinate \in \partial \RefDomain.
\end{align*}
In these experiments a standard VEM of degree $k$ is used to approximate the harmonic extension equations on $\RefDomain$ \cite{daVeiga2013VirtualProblems,Dedner2022ASpaces}. The boundary condition is approximated by considering a piecewise polynomial interpolation of $\partial \Omega_0$.  

The numerical tests are conducted on the same sequence of CVT meshes used in \cite{Wells2024AMethod}, with an initial time step size of $\Delta t = 1e-4$ to a final time of $T = t_0 + 0.01$. The time step size is reduced with each mesh refinement according to $\Delta t \sim h^{k+1}$ to ensure numerical stability. The $L^2$ and $H^1$ error results presented in Table \ref{table::PMEL2H1} suggest that the velocity-based moving mesh method gains a significant improvement in the orders of accuracy for $k=2$ but for $k=3$ is restricted to $O(h^2)$ and $O(h)$ in the $L^2$ and $H^1$ error norms respectively.

\begin{table}[ht]
\centering
\caption{$L^2$ and $H^1$ solution error data for the PME for $k=1,2,3$.}\label{table::PMEL2H1}
\begin{tabular}{|c|c|c|c|c|c|c|}
\hline
\( h \) & \multicolumn{2}{c|}{\( k=1 \)} & \multicolumn{2}{c|}{\( k=2 \)} & \multicolumn{2}{c|}{\( k=3 \)} \\
\cline{2-7}
& L2  & H1 & L2  & H1 & L2  & H1 \\
\hline
0.1661 & 3.254e-3 & 5.899e-2 & 1.622e-3 & 7.689e-2 & 9.911e-4 & 5.470e-2 \\
0.0818 & 9.053e-4 & 2.680e-2 & 2.470e-4 & 2.369e-2 & 2.131e-4 & 2.154e-2 \\
0.0403 & 2.462e-4 & 1.201e-2 & 4.126e-5 & 7.870e-3 & 5.592e-5 & 9.920e-3 \\
\hline
rate   & 1.780 & 1.098 & 2.447 & 1.507 & 1.829 & 1.060 \\
\hline
\end{tabular}
\end{table}

Next, we directly compare the performance of the moving mesh method against the numerical results of the linear VEM studied in \cite{Wells2024AMethod}. The solution and mesh errors are defined at time $t_N = T$ as
\begin{align*}
    \norm{\rho_{h,N} - \rho_N}_{sol} &:= \frac{1}{\Ndofs} \sum_{i=1}^{\Ndofs} \abs{\text{dof}_i(\rho_{h,N}) - \text{dof}_i(\rho_N)}, \\
     \norm{\ALEMap_{h,N} - \ALEMap_N}_{mesh} &:= \frac{1}{N_B} \sum_{i=1}^{N_B} \abs{R_i - r_0 \lambda(t_N)},
\end{align*}
where $\{ R_i \}_{i=1}^{N_B}$ is the set of radial distance from the origin of boundary DoFs of $\ALEMap_{h,N}$ and $\lambda(t)$ is given by Equation \eqref{eq::PMEextras}.

The $l^1$ errors are presented in Table \ref{table::PMEl1} in which the method is seen to achieve $O(h^3)$ accuracy in the solution and $O(h^2)$ accuracy in the mesh for $k=2$.
As with the results of Table \ref{table::PMEL2H1}, we observe the method is limited to $O(h^2)$ and $O(h)$ for the solution and mesh error for $k=3$.
 
\begin{table}[ht]
\centering
\caption{$l^1$ solution and mesh error data for the PME for $k=1,2,3$.}\label{table::PMEl1}
\begin{tabular}{|c|c|c|c|c|c|c|}
\hline
\( h \) & \multicolumn{2}{c|}{\( k=1 \)} & \multicolumn{2}{c|}{\( k=2 \)} & \multicolumn{2}{c|}{\( k=3 \)} \\
\cline{2-7}
& Sol & Mesh & Sol & Mesh & Sol & Mesh \\
\hline
0.1661 & 2.539e-3 & 1.570e-3 & 1.150e-3 & 1.467e-3 & 5.559e-4 & 9.902e-4 \\
0.0818 & 5.976e-4 & 1.069e-3 & 1.276e-4 & 3.873e-4 & 8.332e-5 & 2.834e-4 \\
0.0403 & 1.384e-4 & 5.320e-4 & 1.579e-5 & 1.072e-4 & 1.549e-5 & 1.028e-4 \\
\hline
rate   & 2.000 & 0.954 & 2.857 & 1.757 & 2.300 & 1.387 \\
\hline
\end{tabular}
\end{table}

We conclude this section with a brief discussion on the possible causes of the sub-optimal convergence rates for $k=3$. 
Existing numerical results suggested that sub-optimal convergence rates occurred for $k > 1$ when a simple interpolation of the boundary of $\Omega_0$ are used to define $\ALEMap_{h,0}$ \cite{WellsThesis}. We also note from the numerical experiments of Section \ref{sec::numerics}, where $\Omega_0$ was interpolated exactly, the ALE-VEM achieved the expected $O(h^{k+1})$ and $O(h^k)$ in the $L^2$ and $H^1$ errors respectively. 

To test the impact of the approximation of the initial ALE map $\ALEMap_0$, we consider the velocity-based moving mesh algorithm with weakly imposed boundary conditions and an initial ALE map of the identity
\[
\rho_{h,n} \in V_{h,n}, \qquad \Omega_{h,0} \equiv \RefDomain.
\]
These assumptions leads to a modification of the PME solution which has the same point values given by Equation \eqref{eq::PMEsolution}, but now with a non-zero trace on $\partial \Omega_{h,0}$. As observed in Table \ref{table::PMEweak} this leads to the expected $O(h^{k+1})$ accuracy in the VEM solution for $k=1,2,3$ and an upper limit of $O(h^2)$ for the mesh accuracy due to the linear approximation of $\Omega_0$. 
These results along with experimental results from \cite{WellsThesis} suggest that indeed the choice of approximating the initial domain influences the expected orders of convergence in the ALE-VEM scheme. Possible remedies to this issue could involve solving a different ALE equation for $\ALEMap_{h,0}$ such as bi-harmonic extension \cite{Richter2017}, or by employing a higher regularity VEM such as $C^1$ conforming VEM spaces \cite{Antonietti2016AMeshes}. These methods are beyond the scope of this paper and are left for future investigations.

\begin{table}[ht]
\centering
\caption{$l^1$ solution an mesh error data for the PME with weakly imposed boundary conditions for $k=1,2,3$.}\label{table::PMEweak}
\begin{tabular}{|c|c|c|c|c|c|c|}
\hline
\( h \) & \multicolumn{2}{c|}{\( k=1 \)} & \multicolumn{2}{c|}{\( k=2 \)} & \multicolumn{2}{c|}{\( k=3 \)} \\
\cline{2-7}
& Sol & Mesh & Sol & Mesh & Sol & Mesh \\
\hline
0.1661 & 2.539e-3 & 1.570e-3 & 3.739e-4 & 2.204e-3 & 3.739e-4 & 1.470e-3 \\
0.0818 & 5.976e-4 & 1.069e-3 & 4.931e-5 & 4.389e-4 & 2.464e-5 & 2.926e-4 \\
0.0403 & 1.384e-4 & 5.320e-4 & 6.190e-6 & 1.012e-4 & 1.550e-6 & 6.745e-5 \\
\hline
rate   & 2.066 & 0.986 & 2.931 & 2.072 & 3.907 & 2.073 \\
\hline
\end{tabular}
\end{table}

\section{Conclusion}\label{sec::conclusion}
 In this paper, a high-order virtual element discretisation of a conservative Arbitrary Lagrangian-Eulerian scheme is presented for the convection-diffusion equation. An isoparametric Virtual Element Method was used to formulate a semi-discrete ALE scheme for the convection-diffusion equation. Numerical results presented suggest the method attains optimal high-order convergence rates of $O(h^k)$ and $O(h^{k+1})$ in the $H^1$ and $L^2$ norms respectively. Extending this method to a velocity-based moving mesh algorithm for implicit moving boundary problems leads to improved orders of convergence, attaining optimal $O(h^3)$ accuracy in the solution error using the quadratic VEM.

 Further investigation is required to understand the sub-optimal convergence rates for the implicit moving boundary problems for higher-order VEM discretisations. The development of a three-dimensional ALE-VEM scheme remains an open problem. Advancements in this direction would first require the study of isoparametric VEMs in three dimensions.

\textit{A priori} stability and error estimates for this ALE-VEM scheme are currently being developed. 
Formulating an ALE-VEM scheme for the Incompressible Navier-Stokes equations will be the subject of future studies.

\section*{Acknowledgements}
The author was supported by EPSRC doctoral training grants EP/N50970X/1 and EP/R513283/1.
The author would also like to thank Andrea Cangiani (SISSA), Thomas Radley and Matthew Hubbard (University of Nottingham) for their invaluable discussions and feedback on this work.

\bibliography{References}{}

\begin{thebibliography}{10}

\bibitem{Adak2021AFormulation}
D.~Adak, D.~Mora, S.~Natarajan, and A.~Silgado.
\newblock {A virtual element discretization for the time dependent
  Navier–Stokes equations in stream-function formulation}.
\newblock {\em ESAIM: Mathematical Modelling and Numerical Analysis},
  55(5):2535--2566, 4 2021.

\bibitem{Ahmad2013EquivalentMethods}
B.~Ahmad, A.~Alsaedi, F.~Brezzi, L.~D. Marini, and A.~Russo.
\newblock {Equivalent projectors for virtual element methods}.
\newblock {\em Computers {\&} Mathematics with Applications}, 66(3):376--391, 4
  2013.

\bibitem{Alkamper2016TheModule}
M.~Alk{\"{a}}mper, A.~Dedner, R.~Kl{\"{o}}fkorn, and M.~Nolte.
\newblock {The DUNE-ALUGRID Module}.
\newblock {\em Archive of Numerical Software}, 4(1):1--28, 2016.

\bibitem{Antonietti2016AMeshes}
P.~F. Antonietti, L.~Beir{\~{a}}o Da~Veiga, S.~Scacchi, and M.~Verani.
\newblock {A C1 virtual element method for the Cahn-Hilliard equation with
  polygonal meshes}.
\newblock {\em SIAM Journal on Numerical Analysis}, 54(1):34--56, 4 2016.

\bibitem{Antonietti2014AMeshes}
P.~F. Antonietti, L.~B. da~Veiga, D.~Mora, and M.~Verani.
\newblock {A Stream Virtual Element Formulation of the Stokes Problem on
  Polygonal Meshes}.
\newblock {\em SIAM Journal on Numerical Analysis}, 52(1):386--404, 1 2014.

\bibitem{Antonietti2022VirtualEquation}
P.~F. Antonietti, G.~Vacca, and M.~Verani.
\newblock {Virtual element method for the Navier–Stokes equation coupled with
  the heat equation}.
\newblock {\em IMA Journal of Numerical Analysis}, 2022.

\bibitem{Bachini2021Arbitrary-orderSurfaces}
E.~Bachini, G.~Manzini, and M.~Putti.
\newblock {Arbitrary-order intrinsic virtual element method for elliptic
  equations on surfaces}.
\newblock {\em Calcolo}, 58(3):30, 4 2021.

\bibitem{Baines2005}
M.~Baines, M.~Hubbard, and P.~Jimack.
\newblock {A moving mesh finite element algorithm for the adaptive solution of
  time-dependent partial differential equations with moving boundaries}.
\newblock {\em Applied Numerical Mathematics}, 54(3-4):450--469, 8 2005.

\bibitem{Baines2011}
M.~J. Baines, M.~E. Hubbard, and P.~K. Jimack.
\newblock {Velocity-Based Moving Mesh Methods for Nonlinear Partial
  Differential Equations}.
\newblock {\em Communications in Computational Physics}, 10(3):509--576, 4
  2011.

\bibitem{Baines2006}
M.~J. Baines, M.~E. Hubbard, P.~K. Jimack, and A.~C. Jones.
\newblock {Scale-invariant moving finite elements for nonlinear partial
  differential equations in two dimensions}.
\newblock {\em Applied Numerical Mathematics}, 56(2):230--252, 4 2006.

\bibitem{Baines2009AProblems}
M.~J. Baines, M.~E. Hubbard, P.~K. Jimack, and R.~Mahmood.
\newblock {A moving-mesh finite element method and its application to the
  numerical solution of phase-change problems}.
\newblock {\em Communications in Computational Physics}, 6(3):595--624, 2009.

\bibitem{Bastian2021TheDevelopments}
P.~Bastian, M.~Blatt, A.~Dedner, N.-A. Dreier, C.~Engwer, R.~Fritze,
  C.~Gr{\"{a}}ser, C.~Gr{\"{u}}ninger, D.~Kempf, R.~Kl{\"{o}}fkorn,
  M.~Ohlberger, and O.~Sander.
\newblock {The Dune framework: Basic concepts and recent developments}.
\newblock {\em Computers {\&} Mathematics with Applications}, 81:75--112, 1
  2021.

\bibitem{basicprinciples}
L.~Beir{\~{a}}o~da Veiga, F.~Brezzi, A.~Cangiani, G.~Manzini, L.~D. Marini, and
  A.~Russo.
\newblock {Basic principles of the virtual element method}.
\newblock {\em Mathematical Models and Methods in Applied Sciences},
  23(01):199--214, 1 2013.

\bibitem{DaVeiga2014}
L.~Beir{\~{a}}o~da Veiga, F.~Brezzi, L.~D. Marini, and A.~Russo.
\newblock {The Hitchhiker's Guide to the Virtual Element Method}.
\newblock {\em Mathematical Models and Methods in Applied Sciences},
  24(08):1541--1573, 7 2014.

\bibitem{ellipticVEM}
L.~Beir{\~{a}}o Da~Veiga, F.~Brezzi, L.~D. Marini, and A.~Russo.
\newblock {Virtual Element Method for general second-order elliptic problems on
  polygonal meshes}.
\newblock {\em Mathematical Models and Methods in Applied Sciences},
  26(4):729--750, 2016.

\bibitem{Bonito2013Time-discreteAnalysis}
A.~Bonito, I.~Kyza, and R.~H. Nochetto.
\newblock {Time-discrete higher order ALE formulations: A priori error
  analysis}.
\newblock {\em Numerische Mathematik}, 125(2):225--257, 2013.

\bibitem{Bonito2013Time-DiscreteStability}
A.~Bonito, I.~Kyza, and R.~H. Nochetto.
\newblock {Time-Discrete Higher-Order ALE Formulations: Stability}.
\newblock {\em SIAM Journal on Numerical Analysis}, 51(1):577--604, 4 2013.

\bibitem{Brenner2008TheMethods}
S.~C. Brenner and L.~R. Scott.
\newblock {\em {The Mathematical Theory of Finite Element Methods}}, volume~15
  of {\em Texts in Applied Mathematics}.
\newblock Springer New York, New York, NY, 2008.

\bibitem{Burden2016NumericalEd.}
A.~Burden, R.~L.~Burden, and J.~Douglas~Faires.
\newblock {\em {Numerical Analysis, 10th ed.}}, volume~10.
\newblock 2016.

\bibitem{IsoVEM}
A.~Cangiani, A.~Dedner, M.~Hubbard, and H.~Wells.
\newblock {Isoparametric Virtual Element Methods}.
\newblock {\em arXiv preprint arXiv:2404.11603}, 4 2024.

\bibitem{Cangiani2017Hp-VersionMeshes}
A.~Cangiani, Z.~Dong, E.~H. Georgoulis, and P.~Houston.
\newblock {\em {hp-Version Discontinuous Galerkin Methods on Polygonal and
  Polyhedral Meshes}}.
\newblock Springer Briefs in Mathematics. Springer International Publishing,
  Cham, 2017.

\bibitem{ConNonConVEM}
A.~Cangiani, G.~Manzini, and O.~J. Sutton.
\newblock {Conforming and nonconforming virtual element methods for elliptic
  problems}.
\newblock {\em IMA Journal of Numerical Analysis}, pages 0--36, 8 2016.

\bibitem{daVeiga2013VirtualProblems}
L.~B. da~Veiga, F.~Brezzi, and L.~D. Marini.
\newblock {Virtual Elements for Linear Elasticity Problems}.
\newblock {\em SIAM Journal on Numerical Analysis}, 51(2):794--812, 1 2013.

\bibitem{daVeiga2017DivergenceMeshes}
L.~B. da~Veiga, C.~Lovadina, and G.~Vacca.
\newblock {Divergence free virtual elements for the stokes problem on polygonal
  meshes}.
\newblock {\em ESAIM: Mathematical Modelling and Numerical Analysis},
  51(2):509--535, 3 2017.

\bibitem{daVeiga2018VirtualMeshes}
L.~B. da~Veiga, C.~Lovadina, and G.~Vacca.
\newblock {Virtual Elements for the Navier--Stokes Problem on Polygonal
  Meshes}.
\newblock {\em SIAM Journal on Numerical Analysis}, 56(3):1210--1242, 1 2018.

\bibitem{Dedner2022ASpaces}
A.~Dedner and A.~Hodson.
\newblock {A framework for implementing general virtual element spaces}.
\newblock {\em SIAM Journal of Scientific Computing (to be published)}, 8 2024.

\bibitem{Dedner2020PythonModule}
A.~Dedner, R.~Kl{\"{o}}fkorn, and R.~Kl{\"{o}}fkorn.
\newblock {Python Bindings for the DUNE-FEM module}, 2020.

\bibitem{Dedner2012Dune-Fem:Computing}
A.~Dedner, R.~Kl{\"{o}}fkorn, M.~Nolte, and M.~Ohlberger.
\newblock {Dune-Fem: A General Purpose Discretization Toolbox for Parallel and
  Adaptive Scientific Computing}.
\newblock In {\em Advances in DUNE}, pages 17--31. Springer Berlin Heidelberg,
  Berlin, Heidelberg, 2012.

\bibitem{HHOBook}
D.~A. Di~Pietro and J.~Droniou.
\newblock {\em {The Hybrid High-Order Method for Polytopal Meshes}}, volume~19
  of {\em MS{\&}A}.
\newblock Springer International Publishing, Cham, 2020.

\bibitem{Donea1982AnInteractions}
J.~Donea, S.~Giuliani, and J.~Halleux.
\newblock {An arbitrary lagrangian-eulerian finite element method for transient
  dynamic fluid-structure interactions}.
\newblock {\em Computer Methods in Applied Mechanics and Engineering},
  33(1-3):689--723, 9 1982.

\bibitem{Donea2004}
J.~Donea, A.~Huerta, J.-P. Ponthot, and A.~Rodr{\'{i}}guez-Ferran.
\newblock {Arbitrary Lagrangian-Eulerian Methods}.
\newblock In {\em Encyclopedia of Computational Mechanics}, pages 1--25. John
  Wiley {\&} Sons, Ltd, Chichester, UK, 4 2004.

\bibitem{Du2006ConvergenceTessellations}
Q.~Du, M.~Emelianenko, and L.~Ju.
\newblock {Convergence of the Lloyd Algorithm for Computing Centroidal Voronoi
  Tessellations}.
\newblock {\em SIAM Journal on Numerical Analysis}, 44(1):102--119, 4 2006.

\bibitem{Formaggia1999StabilityElements}
L.~Formaggia and F.~Nobile.
\newblock {Stability analysis for the arbitrary Lagrangian Eulerian formulation
  with finite elements}.
\newblock {\em East-West Journal of Numerical Mathematics}, 7(2):105--131,
  1999.

\bibitem{Formaggia2004StabilityALEFEM}
L.~Formaggia and F.~Nobile.
\newblock {Stability analysis of second-order time accurate schemes for
  ALE–FEM}.
\newblock {\em Computer Methods in Applied Mechanics and Engineering},
  193(39-41):4097--4116, 4 2004.

\bibitem{Gaburro2021AChange}
E.~Gaburro.
\newblock {A Unified Framework for the Solution of Hyperbolic PDE Systems Using
  High Order Direct Arbitrary-Lagrangian–Eulerian Schemes on Moving
  Unstructured Meshes with Topology Change}.
\newblock {\em Archives of Computational Methods in Engineering},
  28(3):1249--1321, 2021.

\bibitem{Gaburro2020HighChanges}
E.~Gaburro, W.~Boscheri, S.~Chiocchetti, C.~Klingenberg, V.~Springel, and
  M.~Dumbser.
\newblock {High order direct Arbitrary-Lagrangian-Eulerian schemes on moving
  Voronoi meshes with topology changes}.
\newblock {\em Journal of Computational Physics}, 407, 2020.

\bibitem{Gastaldi2001AElements}
L.~Gastaldi.
\newblock {A priori error estimates for the Arbitrary Lagrangian Eulerian
  formulation with finite elements}.
\newblock {\em Journal of Numerical Mathematics}, 9(2):123--156, 1 2001.

\bibitem{Hirt1974AnSpeeds}
C.~W. Hirt, A.~A. Amsden, and J.~L. Cook.
\newblock {An arbitrary Lagrangian-Eulerian computing method for all flow
  speeds}.
\newblock {\em Journal of Computational Physics}, 14(3), 1974.

\bibitem{MMFEMBoundaries}
M.~E. Hubbard, M.~J. Baines, and P.~K. Jimack.
\newblock {Consistent Dirichlet boundary conditions for numerical solution of
  moving boundary problems}.
\newblock {\em Applied Numerical Mathematics}, 59(6), 2009.

\bibitem{Irisarri2019StabilizedEquations}
D.~Irisarri and G.~Hauke.
\newblock {Stabilized virtual element methods for the unsteady incompressible
  Navier–Stokes equations}.
\newblock {\em Calcolo}, 56(4), 2019.

\bibitem{lipnikov2D}
K.~Lipnikov and N.~Morgan.
\newblock {A high-order conservative remap for discontinuous Galerkin schemes
  on curvilinear polygonal meshes}.
\newblock {\em Journal of Computational Physics}, 399, 2019.

\bibitem{lipnikov2020ConservativeMeshes}
K.~Lipnikov and N.~Morgan.
\newblock {Conservative high-order discontinuous Galerkin remap scheme on
  curvilinear polyhedral meshes}.
\newblock {\em Journal of Computational Physics}, 420:109712, 4 2020.

\bibitem{Liu2019TheEquations}
X.~Liu and Z.~Chen.
\newblock {The nonconforming virtual element method for the Navier-Stokes
  equations}.
\newblock {\em Advances in Computational Mathematics}, 45(1), 2019.

\bibitem{Marlow2011MovingEquations}
R.~Marlow, M.~E. Hubbard, and P.~K. Jimack.
\newblock {Moving mesh methods for solving parabolic partial differential
  equations}.
\newblock {\em Computers {\&} Fluids}, 46(1):353--361, 4 2011.

\bibitem{Mazzia2020VirtualEnvironments}
A.~Mazzia, M.~Ferronato, P.~Teatini, and C.~Zoccarato.
\newblock {Virtual element method for the numerical simulation of long-term
  dynamics of transitional environments}.
\newblock {\em Journal of Computational Physics}, 407, 2020.

\bibitem{Ngo2017AEquation}
C.~Ngo and W.~Huang.
\newblock {A study on moving mesh finite element solution of the porous medium
  equation}.
\newblock {\em Journal of Computational Physics}, 331:357--380, 4 2017.

\bibitem{Richter2017}
T.~Richter.
\newblock {\em {Fluid-structure Interactions}}, volume 118 of {\em Lecture
  Notes in Computational Science and Engineering}.
\newblock Springer International Publishing, Cham, 2017.

\bibitem{Senechal1995SpatialDiagrams}
M.~Senechal, A.~Okabe, B.~Boots, and K.~Sugihara.
\newblock {Spatial Tessellations: Concepts and Applications of Voronoi
  Diagrams}.
\newblock {\em The College Mathematics Journal}, 26(1):79, 1995.

\bibitem{Souli2000ALEProblems}
M.~Souli, A.~Ouahsine, and L.~Lewin.
\newblock {ALE formulation for fluid-structure interaction problems}.
\newblock {\em Computer Methods in Applied Mechanics and Engineering},
  190(5-7), 2000.

\bibitem{Sutton2017TheMATLAB}
O.~J. Sutton.
\newblock {The virtual element method in 50 lines of MATLAB}.
\newblock {\em Numerical Algorithms}, 75(4):1141--1159, 4 2017.

\bibitem{Takashi1992AnBody}
N.~Takashi and T.~J. Hughes.
\newblock {An arbitrary Lagrangian-Eulerian finite element method for
  interaction of fluid and a rigid body}.
\newblock {\em Computer Methods in Applied Mechanics and Engineering},
  95(1):115--138, 2 1992.

\bibitem{Talischi2012PolyMesher:Matlab}
C.~Talischi, G.~H. Paulino, A.~Pereira, and I.~F.~M. Menezes.
\newblock {PolyMesher: a general-purpose mesh generator for polygonal elements
  written in Matlab}.
\newblock {\em Structural and Multidisciplinary Optimization}, 45(3):309--328,
  4 2012.

\bibitem{Vacca2017VirtualMeshes}
G.~Vacca.
\newblock {Virtual Element Methods for hyperbolic problems on polygonal
  meshes}.
\newblock {\em Computers and Mathematics with Applications}, 74(5):882--898,
  2017.

\bibitem{Vacca2015VirtualMeshes}
G.~Vacca and L.~Beir{\~{a}}o~da Veiga.
\newblock {Virtual element methods for parabolic problems on polygonal meshes}.
\newblock {\em Numerical Methods for Partial Differential Equations},
  31(6):2110--2134, 11 2015.

\bibitem{Vazquez2006TheEquation}
J.~L. Vazquez.
\newblock {\em {The Porous Medium Equation}}.
\newblock Oxford University Press, 4 2006.

\bibitem{Verma2023VirtualMeshes}
N.~Verma and S.~Kumar.
\newblock {Virtual element approximations for non-stationary navier-stokes
  equations on polygonal meshes}.
\newblock {\em Journal of Applied Analysis and Computation}, 13(3), 2023.

\bibitem{Wang2019AMeshes}
G.~Wang, F.~Wang, L.~Chen, and Y.~He.
\newblock {A divergence free weak virtual element method for the Stokes–Darcy
  problem on general meshes}.
\newblock {\em Computer Methods in Applied Mechanics and Engineering},
  344:998--1020, 2 2019.

\bibitem{WellsThesis}
H.~Wells.
\newblock {\em {Moving Mesh Virtual Element Methods}}.
\newblock PhD thesis, University of Nottingham, 2023.

\bibitem{Wells2024AMethod}
H.~Wells, M.~E. Hubbard, and A.~Cangiani.
\newblock {A velocity-based moving mesh virtual element method}.
\newblock {\em Computers {\&} Mathematics with Applications}, 155:110--125, 4
  2024.

\bibitem{Wesseling2001PrinciplesDynamics}
P.~Wesseling.
\newblock {\em {Principles of Computational Fluid Dynamics}}, volume~29 of {\em
  Springer Series in Computational Mathematics}.
\newblock Springer Berlin Heidelberg, Berlin, Heidelberg, 2001.

\bibitem{Zhao2020TheProblem}
J.~Zhao, B.~Zhang, S.~Mao, and S.~Chen.
\newblock {The nonconforming virtual element method for the Darcy-Stokes
  problem}.
\newblock {\em Computer Methods in Applied Mechanics and Engineering}, 370, 4
  2020.

\end{thebibliography}
\bibliographystyle{abbrv}

\end{document}